%% file: scriptie2.tex
\documentclass[a4paper]{article}

\usepackage{latexsym}
\usepackage{amsmath}
\usepackage{amsthm}
\usepackage{amssymb}
\usepackage{psfig,float}
\usepackage[line,matrix,arrow]{xy}

\theoremstyle{plain}
  \newtheorem{stelling}{Theorem}[subsection]
  \newtheorem{prop}[stelling]{Proposition}
  \newtheorem{lemma}[stelling]{Lemma}
 \newtheorem{gevolg}[stelling]{Corollary}
\theoremstyle{definition} 
  \newtheorem{definitie}[stelling]{Definition}
  \newtheorem{voorbeeld}[stelling]{Example}
 
  \newcommand\textdef{\emph}
\theoremstyle{remark}
  \newtheorem{opmerking}[stelling]{Remark}

\theoremstyle{plain}

\theoremstyle{definition}

\theoremstyle{remark}

\newenvironment{bewijs}[1][Proof]%
  {\begin{proof}[{\upshape\bfseries#1}]}%
  {\end{proof}}

\numberwithin{equation}{section}

\newcommand\m{\mathfrak{M}}
\newcommand\n{\mathfrak{N}}
\newcommand\p{\mathfrak{P}}
\newcommand\q{\mathfrak{Q}}
\newcommand\sqtimes{\boxtimes}
\newcommand\biso{\stackrel{b}{\simeq}}
\newcommand\meq{\stackrel{M}{\sim}}

\begin{document}

\title{{A Bicategorical Approach to Morita Equivalence for Rings and 
von Neumann Algebras}\thanks{This paper is the master's thesis I wrote under supervision of N.P. Landsman, University of Amsterdam. An abbreviated version will appear in the Journal of Mathematical Physics.} }
\author{R.M. Brouwer\thanks{Universiteit van Amsterdam, Korteweg-de Vries Instituut,
Faculteit Natuurwetenschappen, Wiskunde en Informatica,
Plantage Muidergracht 24, 1018 TV Amsterdam}
\thanks{Current address: CWI, Postbus 94079, 1090 GB, Amsterdam. Email: rachel.brouwer@cwi.nl}}

\maketitle


\section*{Abstract}
Rings form a bicategory $\mathsf{[Rings]}$, with classes of bimodules as horizontal arrows, and bimodule maps as vertical arrows. The notion of Morita equivalence for rings can be translated in terms of bicategories in the following way. Two rings are Morita equivalent if and only if they are isomorphic objects in the bicategory. 
We repeat this construction for von Neumann algebras. Von Neumann algebras form a bicategory $\mathsf{[W^*]}$, with classes of correspondences as horizontal arrows, and intertwiners as vertical arrows. 
Two von Neumann algebras are Morita equivalent if and only if they are isomorphic objects in the bicategory $\mathsf{[W^*]}$. 
\count0=1
\section*{Introduction}
\addcontentsline{toc}{section}{Introduction}

This paper is concerned with the 
 Morita theory of rings and von Neumann algebras. 
Abstract ring theory was initiated around $1920$, by, amongst others, 
Fraenkel, Brauer, Artin, Hasse, and Emmy Noether. 
Following from abstract field theory, ring theory has found its use in many areas of mathematics. Von Neumann algebras, first introduced by von Neumann in $1930$, are now widely used in analysis and mathematical physics.

At first sight, the rather abstract field of ring and module theory, and the more physical von Neumann algebra theory are wide apart. For example, ring theory is used in number theory and algebraic geometry. Von Neumann algebras find their use in quantum field theory and statistical mechanics, as well as in representation theory and ergodic theory. In this paper,
 the author has tried to show analogies between the two, using Morita equivalence and bicategories. However, bicategories, though interesting objects as such, are nevertheless just a tool for handling Morita equivalence.

Kiiti Morita introduced the traditional notion of Morita equivalence for rings, where it is said that two rings $R,S$ are Morita equivalent if 
their categories of right modules, $\mathfrak{M}_R$ and $\mathfrak{M}_S$, are equivalent. 
In the case of von Neumann algebras, $\m,\n$ are said 
to be Morita equivalent if there exists a correspondence $\m  \rightarrow \mathcal{H} \leftarrow \n$ for which the representation of $\m$ on 
$\mathcal{H}$ is faithful and \\$\m \simeq (\n^{op})^{\prime}$ holds. 
However, we may choose a definition of Morita equivalence for von Neumann algebras, similar to the definition of Morita equivalence for rings (in terms of 
representation categories), which is equivalent to our definition.

The main result of this paper is the following:
For both rings and von Neumann algebras, it is possible 
to prove that Morita equivalence is nothing but isomorphism in their respective bicategories. Despite the fact that rings and von Neumann algebras have their use in different areas of mathematics, they have the same underlying structure as far as Morita equivalence is concerned. 
\ \\ \ \\
This paper is organized as follows. 
Section $1$ contains the basic definition of a bicategory. A few examples will be discussed. 

Section $2$ handles the case of rings. First, the tensor product of two modules will be explained in Subsection $2.1$. Next, 
we will consider the traditional Morita theory, which handles progenerators and categories of modules in Subsection $2.2$. 
The notion of a bicategory emerges in Subsection $2.3$, where we will show that rings form a bicategory, and we will state the 
Morita theory in terms of bicategories. To justify the use of the terminology ``Morita theory'' in the bicategory case, we will show that the 
results of Subsections $2.2$ and $2.3$ lead to equivalent theories. This is done in Subsection $2.4$.   

Section $3$ handles von Neumann algebras. The goal of this section is to explain the notion 
of Morita equivalence, and to exhibit corresponding statements in terms of bicategories. It turns out to be the case 
that for the construction of the bicategory of von Neumann algebras, one needs the standard form and the identity 
correspondence of a von Neumann algebra, as well as the concept of Connes fusion. These are discussed in Subsections 
$3.1,3.2$ and $3.3$. Finally, one will find the Morita theory of von Neumann algebras in Subsection $3.4$.
\ 
\include{bicategories}
\include{rings}
\include{vNa}

\include{bibl}
\end{document}

%% file: bicategories.tex
\begin{section}{Bicategories}
\begin{subsection}{Definitions}
This section will explain the notion of a \emph{bicategory} 
that we will use in later sections. It is assumed that the reader is familiar with categories, functors and natural equivalences or natural isomorphisms. 
In this paper, we will use the convention that the class of objects of a 
category ${C}$ is denoted by $C_0$; the 
class of morphisms of $C$ is denoted by $C_1$. The notation $(A,B)$ is used for all arrows $B \rightarrow A$, which allows us to write the composition of arrows conveniently.
For a standard text on categories the reader is referred to~\cite{maclane}. 
See also~\cite{ben} and~\cite{maclane} for an overview of bicategories.

In several situations where we have a bifunctor $B \times B \rightarrow B$ 
(for a category $B$), this bifunctor is not associative. 
If it is, and has a unit element, our category $B$ becomes a so-called strict monoidal category. For a 
(relaxed) monoidal category, there exist natural equivalences such that the 
bifunctor is associative up to isomorphism. Further, sometimes we would like 
to define bifunctors $B \times C \rightarrow D$ for categories $B,C,D$. Such 
bifunctors or composition functors give rise to a $2$-category if they are associative. 
Again, composition functors are generally not associative. An example is a bicategory, where the composition functor is merely associative up to isomorphism.

\begin{definitie}
A \textdef{bicategory} $\mathcal{B}$ consists of the following ingredients: 
\begin{enumerate}
\item{A set $\mathcal{B}_{0}$ of objects.}
\item{For all pairs $(A,B)$ of objects, a category. If there is no confusion 
possible, this category will also be denoted by $(A,B)$. The class of all such categories will be denoted by $\mathcal{B}_{1}$; 
it contains all morphisms or \textdef{horizontal arrows} of $\mathcal{B}$. 
The morphisms (arrows) $(A,B)_1$ of the category $(A,B)$ are called \textdef{vertical arrows}.} 
\item{For each triple $(A,B,C)$ of objects in  $B_{0}$, 
a composition functor.

\begin{equation}
 \mathcal{C}(A,B,C) : (A,B) \times (B,C) \longrightarrow (A,C). 
\end{equation}
If $(P,Q)$ is an element of $(A,B) \times (B,C)$, we will write $P \ast Q$ for\\ $\mathcal{C}(A,B,C)(P,Q)$. The same notation will be used on the arrows.
}
\item{For each object $B$ of $\mathcal{B}_{0}$, an object $I_{B}$ of $(B,B)$. 
$I_{B}$ is called the \textdef{identity arrow} of $B$.}
\item{For each quadruple $(A,B,C,D)$ of objects in $\mathcal{B}_{0}$, a natural isomorphism $\beta$ between the functors 
\begin{equation}
F = \bigl[\mathcal{C}(A,B,D)\bigr] \circ \bigl[Id_{(A,B)} \times \mathcal{C}(B,C,D)\bigr], 
\end{equation}
and 
\begin{equation}
G = \bigl[\mathcal{C}(A,C,D)\bigr] \circ \bigl[\mathcal{C}(A,B,C) \times Id_{(C,D)}\bigr], 
\end{equation} 
where 
\begin{equation}
 F, G : (A,B) \times (B,C) \times (C,D) \longrightarrow (A,D). 
\end{equation}
This means there exists a natural equivalence $\beta$ which assigns to every object $(E_1,E_2,E_3)$ in  $(A,B) \times (B,C) \times (C,D)$ an arrow 
$$\beta(E_1,E_2,E_3) : F_0(E_1,E_2,E_3) \rightarrow G_0(E_1,E_2,E_3),$$ such that
for every arrow $$h : (E_1,E_2,E_3) \rightarrow ({E_1}^{\prime},{E_2}^{\prime},{E_3}^{\prime}) \in 
(A,B) \times (B,C) \times (C,D),$$ the diagram 
\begin{equation}
\xymatrix{
F_0(E_1,E_2,E_3)\ar[rr]^{\beta(E_1,E_2,E_3)}\ar[d]_{F_1(h)} & & G_0(E_1,E_2,E_3)\ar[d]^{G_1(h)} \\
F_0({E_1}^{\prime},{E_2}^{\prime},{E_3}^{\prime}) \ar[rr]^{\beta({E_1}^{\prime},{E_2}^{\prime},{E_3}^{\prime})} & & G_0({E_1}^{\prime},{E_2}^{\prime},{E_3}^{\prime})\\
} 
\end{equation}
commutes. Note that $\beta$ depends on the quadruple $(A,B,C,D)$, but when no confusion is possible, we will omit this.
} 
\item{For each pair of objects $(A,B)$, two natural equivalences $R_{(A,B)}$ 
and $L_{(A,B)}$ called \textdef{left} and \textdef{right identities}. Here $L_{(A,B)}$ is a natural equivalence  between 
\begin{equation}
 \bigl[\mathcal{C}(A,A,B)\bigr] \circ \bigl[I_A \times Id_{(A,B)} \bigr],
\end{equation}
and the canonical functor from $ 1 \times (A,B)$ to $(A,B)$. $R_{(A,B)}$ is a natural equivalence between 
\begin{equation}
\bigl[\mathcal{C}(A,B,B)\bigr] \circ \bigl[Id_{(A,B)} \times I_B \bigr],
\end{equation}
 and the canonical functor from $ (A,B) \times 1$ to $(A,B)$. When no confusion is possible, the subscript of the right and left identities will be omitted. 
}
\item{The natural isomorphisms $\beta$, $L_{(A,B)}$ and $R_{(A,B)}$ are required to satisfy the so-called \emph{coherence axioms}.
\begin{itemize}
\item{\emph{Associativity coherence}:\\
If $(P,Q,R,S)$ is an object in $(A,B) \times (B,C) \times (C,D) \times (D,E)$, the following diagram commutes:
\begin{equation}
\xymatrix{
((P\ast Q)\ast R)\ast S \ar[rr]^{\beta(P,Q,R) \ast Id} \ar[d]_{\beta(P \ast Q, R,S)} & & (P \ast (Q \ast R)) \ast S \ar[d]^{\beta(P,Q \ast R,S)} \\
(P \ast Q) \ast (R \ast S) \ar[dr]_{\beta(P,Q,R\ast S)\ \ \ } & & P \ast (( Q \ast R) \ast S) \ar[dl]^{\ \ \ Id \ast \beta(Q,R,S)} \\
& P \ast ( Q \ast (R \ast S)). & \\
}
\end{equation} 
}
\item{\emph{Identity coherence}:\\
If $(P,Q)$ is an object in $(A,B) \times (B,C)$ the following diagram commutes:
\begin{equation}
\xymatrix{
(P \ast I_B) \ast Q \ar[rr]^{\beta(P,I_B,Q)}\ar[dr]_{R(P)\ast Id\ \ }  &  & P\ast (I_B \ast Q)\ar[dl]^{\ \ Id \ast L(Q)} \\
& P \ast Q. & \\
}
\end{equation}
}
\end{itemize}
}
\end{enumerate}
\end{definitie}
One would like that all diagrams concerning $\beta, L$ and $R$ are commutative.
In fact, if the coherence axioms are satisfied, all such diagrams commute. 
(\!~\cite{maclane} Section XI.7; c.f. the coherence axioms in the case of a monoidal category,~\cite{maclane}, Section VII.2, where commutativity of all diagrams is proven.)

Considering a bicategory and taking isomorphism classes in the categories $(A,B)$ as arrows, we get a category. The natural isomorphisms $L$ and $R$ provide 
left and right identities, the natural isomorphism $\beta$ provides associativity. The coherence axioms even make this category monoidal.
A \hbox{$2$-category}, where the composition is associative, 
is a special case of a bicategory: In this case 
the natural isomorphisms $\beta, L$ and $R$ are all identities, so the coherence axioms are satisfied immediately. 
\label{bic}

For later use, we need a notion of isomorphism in a bicategory, which is 
broader than the usual notion of isomorphism in category theory.. 
\begin{definitie}
Two objects $A,B$ in a bicategory $\mathcal{B}$ are \textdef{isomorphic in the 
bicategory}, $A \biso B$, if an invertible (horizontal) arrow $P$ exists, $P \in (A,B)$, 
i.e.  
\begin{equation}
P \ast P^{-1} \cong I_A,
\end{equation}
and
\begin{equation}
P^{-1} \ast P \cong I_B,
\end{equation}
where the symbol $\cong$ denotes isomorphism of objects in $(A,A)$ and in $(B,B)$ respectively, in the usual sense of categories.
\label{defisobic}
\end{definitie}

Note that every object of the bicategory is isomorphic to itself via its
identity arrow, which is its own inverse. Further, two objects $A,B$ that are isomorphic in the usual sense, are also 
isomorphic in the bicategory, since in that case the categories $(A,A)$ and $(A,B)$ are equivalent. Hence we have isomorphism on objects through the natural 
equivalence. 
The invertible arrow in this case is given by the image of $I_A$ in $(A,B)$.
\end{subsection}
\begin{subsection}{Examples}
\begin{itemize}
\item{
All $2$-categories are bicategories. For example:
\begin{itemize}
\item{The class of all categories as objects, functors as horizontal arrows, and natural transformations as vertical arrows.}
\item{The class of all topological spaces as objects, continuous maps as horizontal arrows, and homotopy classes of continuous maps as vertical arrows.}
\end{itemize}
}
\item{
Each (relaxed) monoidal category $M$ forms a bicategory, which in general is not a \hbox{$2$-category}. The bicategory consists of one object (namely $M$), the 
objects of the category $M$ form the horizontal arrows of the bicategory. A composition functor $M \times M \rightarrow M$ which is associative 
up to isomorphism exists, since $M$ is monoidal. The arrows $M_1$ form the vertical arrows of the bicategory. The natural isomorphisms that 
are associated to the monoidal category make that the coherence axioms are satisfied.
}
\item{
More instructive examples of bicategories are the bicategory $\mathsf{[Rings]}$ and the bicategory $\mathsf{[W^*]}$. The bicategory  $\mathsf{[Rings]}$ 
consists of rings as objects, categories of bimodules as horizontal arrows, and linear maps as vertical arrows. The bicategory $\mathsf{[W^*]}$ consists of 
von Neumann algebras as objects, categories of correspondences as horizontal arrows, and intertwiners as vertical arrows. However, it is not easy to show that 
$\mathsf{[Rings]}$ and  $\mathsf{[W^*]}$ are indeed bicategories. Especially the definition of the composition functor is not trivial. Therefore a proof can be found 
in later sections (see Propositions \ref{rings=bic} and \ref{vna=bic}). However, these examples are the main reason to discuss bicategories, since we use bicategories to show that the notions of Morita equivalence for rings and von Neumann have the same underlying structure.
}
\end{itemize}

\end{subsection}

\end{section}

%% file: rings.tex
\begin{section}{Morita theory for rings}
It is assumed that the reader is familiar with the notion of modules and bimodules of rings. If not, basic ring and module theory may be found in~\cite{cohn2}. 
Throughout this section, all rings will have a unit. Let $R$, $S$ be rings. 
A left $R$-module $M$ will be denoted by $_{R}M$. Right modules will be denoted analogously by $M_R$. A $R$-$S$ bimodule $N$ 
will be denoted by $_{R}N_{S}$ or $R \rightarrow N \leftarrow S$. Both notations will be used in the following.

In this section, we will present the "traditional" Morita theory.
After that, this theory will be reformulated in terms of bicategories. 
Finally, it will be shown that these theories are equivalent. 
In both approaches, the notion of a tensor product of two (bi)modules is needed, so first of all,
we will discuss the tensor product. See~\cite{cap} for an extensive discussion of the tensor product.
\begin{subsection}{The tensor product of bimodules}

To define the tensor product of two modules, we need the following definitions. 
Let $M$, $N$ and $L$ be abelian groups, $R$, $S$ and $T$ rings. 
A map $\psi : M \times N \rightarrow L$ is called \emph{bilinear} if it satisfies $\psi(m + m^{\prime} , n) = \psi(m , n) + \psi(m^{\prime} , n)$ 
and $\psi(m , n + n^{\prime}) = \psi(m , n) + \psi(m , n^{\prime})$, for $m,m^{\prime} \in M, n,n^{\prime} \in N$. 
If $M$ is a left $R$-module, $N$ a right $T$-module and $L$ a $R$-$T$ bimodule, a linear map $\psi : M \times N \rightarrow L$ is 
called \emph{$R$-$T$ linear} if the map $\psi$ intertwines the $R$ and $T$ actions. Further, if $M$ is a right $S$-module and $N$ is a 
left $S$-module, a map $\psi : M \times N \rightarrow L$ is called \emph{$S$-balanced} if it satisfies $\psi(m,sn) = \psi(ms,n)$, for $m \in M, n \in N, s \in S$. 
\begin{prop}
Given three rings $R$, $S$ and $T$, and two bimodules $_{R}M_{N}$ and $_{S}N_{T}$, there exists an $R$-$T$ bimodule 
$_{R}(M \otimes_S N)_T$ and an $S$-balanced $R$-$T$ bilinear map $\tau : M \times N \rightarrow\  _{R}(M \otimes_S N)_T$ 
with the following universal property: For every $R$-$T$ bimodule $L$ and every $S$-balanced $R$-$T$ bilinear map $\phi : M \times N \rightarrow L$ there exists a unique 
$R$-$T$ bilinear map $\alpha : M \otimes_S N \rightarrow L$ such that $\phi = \alpha \circ \tau$. In a commutative diagram:
\begin{equation}
\xymatrix{
M \times N \ar[rr]^{\phi}\ar[rd]^{\tau} & & L \\
& M \otimes_S N. \ar[ur]^{!\alpha} & \\
}
\end{equation} 
\begin{bewijs}

Existence of $M\otimes_S N$ follows directly by construction. Consider $Y$, the free $\mathbb{Z}$-module on $M\times N$ with embedding 
$i : M\times N \rightarrow Y$. Now quotient to $I$, which is generated by elements of the form $(m + m^{\prime}n) - (m,n) - (m^{\prime},n)$ or $(m, n+ n^{\prime}) -
(m ,n ) - (m,n^{\prime})$ or $(ms,n) - (m,sn)$, for $ m, m^{\prime} \in M,\  n, n^{\prime} \in N,\ s \in S$. 
Let $\pi : Y \rightarrow Y/I$ be the canonical surjection. The remaining quotient $Y/I$ forms the bimodule tensor product
 $M\otimes_SN$ of $_{R}M_S$ and $_{S}N_T$. One has to show that this tensor product is an $R$-$T$ bimodule. The left $R$-action on $M$ is defined on 
$M \otimes_SN$ by $r(m\otimes_Sn) := (rm \otimes_Sn)$, which is defined because of the left $R$-action on $M$. Note that $I$ is closed under the left $R$-action, so the $R$-action passes to the quotient. The right $T$-action descends to the quotient likewise. 
Define $\tau = \pi \circ i$. The structure of $M\otimes_SN$ causes $\tau$ to be $S$-balanced and $R$-$T$ bilinear.

Let $L$ be a $R$-$T$ bimodule and let $\phi : M\times N \rightarrow L$ be an $S$-balanced, $R$-$T$ linear map. We now obtain the following diagram:
\begin{equation}
\xymatrix{
M\times N \ar[r]^{\phi}\ar[d]^{i} & L \\
Y\ar[r]^{\pi}\ar@{.>}[ur]^{A} & M \otimes_S N.\ar@{.>} [u]^{\alpha} \\
}
\end{equation}
The universal property of a free module provides the decomposition $\phi = A \circ i$, where $A$ is a linear map. Further, $A$ vanishes on the elements of $I$,
 since $\phi$ is bilinear and $S$-balanced, so we have a decomposition $A = \alpha \circ \pi$, where $\alpha$ is a linear map. 
Trivially, $\alpha$ is bilinear and $S$-balanced because of the definition of $M\otimes_SN$. We need to show that $\alpha$ is an $R$-$T$ bilinear map:
\begin{eqnarray}
\alpha(r(m\otimes_Sn)t) &=& \alpha(rm\otimes_Snt) = \alpha(\tau(rm,nt)) \nonumber \\ &=& \phi(rm,nt) = r\phi(m,n)t 
\nonumber \\ &=& r(\alpha(m\otimes_Sn))t,
\end{eqnarray} for 
$m \in M, n \in N ,r \in R, t\in T$.

Uniqueness of $\alpha$ follows from taking a second $R$-$T$ bilinear $S$-balanced map $\tilde{\alpha}$ such that 
$\phi = \tilde{\alpha} \circ \tau$. Defining $\tilde{A} = \tilde{\alpha} \circ \pi$, a similar diagram as above can be formed,
 with $\tilde{\alpha}$ instead of $\alpha$ and $\tilde{A}$ instead of $A$. However, the universal property of a free module guarantees that $A$ is unique, 
so $\tilde{A} = A$. Now $\tilde{\alpha} \circ \pi = \alpha \circ \pi$ and surjectivity of $\pi$ shows $\tilde{\alpha} = \alpha$.

Finally, we show uniqueness of the pair $(M\otimes_S N,\tau)$. Suppose we have a second pair $(\widetilde{M\otimes_SN},\tilde{\tau})$ 
that satisfies the properties stated in the proposition above. Applying the 
universal property to both $(M\otimes_S N,\tau)$ and $(\widetilde{M\otimes_SN},\tilde{\tau})$, we obtain the following commutative diagram, where both 
$\alpha$ and $\tilde{\alpha}$ are unique:
\begin{equation}
\xymatrix{
& M \otimes_S N \ar[dd]^{!\alpha}\ar[dd] \\
M \times N \ar[ur]^{\tau} \ar[dr]^{\tilde{\tau}} \\
& \widetilde{M \otimes_S N}. \ar[uu]\ar[uu]^{!\tilde{\alpha}}\\
}
\end{equation}
We obtain
\begin{equation}
\left.
  \tau = \tilde{\alpha} \circ \tilde{\tau} \atop {
  \tilde{\tau} = \alpha \circ \tau}
\right\}\Rightarrow
 { \tau = \tilde{\alpha} \circ \alpha \circ \tau; \label{tauid} \atop{
  \tilde{\tau} = \alpha \circ \tilde{\alpha} \circ \tilde{\tau}.}}
\end{equation}
By construction, $\tau$ is surjective, so equation \eqref{tauid} implies $\tilde{\alpha}\circ \alpha = id$ on $M \otimes_S N$. 
Hence $\alpha \circ \tilde{\alpha} = id$ on $\widetilde{M \otimes_S N}$. 
Now $\tilde{\alpha} = \alpha^{-1}$ so $M \otimes_S N \cong \widetilde{M \otimes_S N}$ and $\tau = \alpha^{-1}\circ \tilde{\tau}$.

\end{bewijs}

\end{prop}

\begin{opmerking}
Note that, given the above proposition, we are able to construct the tensor product $M \otimes_R N$  between a left $R$-module $N$ and a 
right $R$-module $M$, for any ring $R$. By considering the left $R$-module $N$ as a $R$-$\mathbb{Z}$ bimodule and the  
right $R$-module $M$ as a $\mathbb{Z}$-$R$ bimodule, the proof of the above proposition applies.
\label{simpeltens}
\end{opmerking}

Let $R,S,T$ be rings. 
As a preparation for the categorical statements of Morita theory,
 we will show that the tensor product $\otimes_S$ defines a functor. 
Let $(R,S)$ denote the following category:
The class of objects $(R,S)_0$ consists of $R$-$S$ bimodules, the class of arrows $(R,S)_1$ consists of $R$-$S$ linear maps. 
The categories $(S,T)$ and $(R,T)$ are defined likewise.
  
On objects, $\otimes_S$ is defined by: 
\begin{eqnarray} 
 \otimes_S : (R,S) \times (S,T) &\longrightarrow& (R,T) \nonumber \\
  _{R}M_S \times _{S}N_T &\longmapsto& _{R}(M\otimes_S N)_T,
\end{eqnarray} 
for $_{R}M_S \in (R,S)$, $ _{S}N_T \in (S,T)$.\\
 On arrows, $\otimes_S$ acts as follows: Let $\phi :{}_{R}M_S \rightarrow {}_{R}K_S$ be an arrow in $(R,S)$, and 
$\psi : {}_SN_T \rightarrow {}_SL_T$ be an arrow in $(S,T)$. Then, for \hbox{$\sum_i(m_i \otimes_S n_i)$} $\in M\otimes_S N$:
\begin{eqnarray}
(\phi \otimes_S \psi) : _{R}(M \otimes_S N)_T &\longrightarrow& _{R}(K \otimes_S L)_T \nonumber \\  
\sum_i(m_i\otimes_S n_i) &\longmapsto& \sum_i(\phi(m_i) \otimes_S \psi(n_i)).
\end{eqnarray}
Since both $\phi$ and $\psi$ intertwine  the $S$-action, $\phi \otimes_S \psi$ is well-defined.

The construction of the tensor product (Proposition \ref{er=tens}) shows that on objects, the image of $\otimes_S$ lies in $(R,T)_0$. On 
arrows, one has to show that for each pair of arrows $(\phi \times \psi) \in \bigl((R,S) \times (S,T)\bigr)_1$, the image $(\phi \otimes_S \psi)$ is an $R$-$T$ linear map.
For $r \in R$, and  $\sum_i (m_i \otimes_S n_i) \in M\otimes_S N$, one has 
\begin{eqnarray}
r(\phi \otimes_S \psi) \sum_i (m_i \otimes_S n_i) &=& r \sum_i \bigl(\phi(m_i) \otimes_S \psi(n_i)\bigr) \nonumber \\
&=&\sum_i \bigl(r\phi(m_i) \otimes_S \psi(n_i)\bigr) \nonumber \\
&=&\sum_i \bigl(\phi(rm_i) \otimes_S \psi(n_i)\bigr) \nonumber \\
&=&(\phi \otimes_S \psi)\sum_i(rm_i \otimes_S n_i) \nonumber \\
&=&(\phi \otimes_S \psi)\Bigl(r\sum_i(m_i \otimes_S n_i)\Bigr).
\end{eqnarray}
A similar computation shows that $(\phi \otimes_S \psi)$ preserves the right $T$-action. 

It is left to show that $\otimes_S$ is a functor. By definition, 
a functor $F$ from a category $\mathcal{C} \rightarrow \mathcal{D}$ is a map which assigns to each object in $\mathcal{C}$ an object in 
$\mathcal{D}$ and to each arrow $f: c \rightarrow c'$ in $\mathcal{C}$ an arrow $F(f) : F(c) \rightarrow F(c')$ in $\mathcal{C}$ such that 
\begin{equation}
F_{id_c} = id_{F(c)},
\label{propfunct1}
\end{equation}
and
\begin{equation}
F(f)\circ F(g) = F(f\circ g),
\label{propfunct2}
\end{equation}
for all objects $c$ in $\mathcal{C}_0$ and all arrows  $f,g$ in $\mathcal{C}_1$, whenever the composition of arrows $f\circ g$ is defined in $\mathcal{C}_1$. 
Let $_{R}M_S \in (R,S)$ and $ _SN_T \in (S,T)$. Then
\begin{eqnarray}
\otimes_S \circ id_{(R,S) \times (S,T)}\bigl(M,N\bigr) &=& \otimes_S 
\bigl(M,N\bigr) \nonumber \\
&=&  _{R}(M\otimes_S N)_T \nonumber \\
&=& id_{(R,T)}(M \otimes_S N) \nonumber \\
&=& id_{\otimes_S\bigl((R,S) \times (S,T)\bigr)}\bigl(M,N\bigr),\nonumber \\
&&
\end{eqnarray}
and
\begin{eqnarray}
\otimes_S ((f_1 \times f_2) \circ (g_1 \times g_2))\sum_i(m_i \otimes_S n_i) &=& ((f_1 \circ g_1) \otimes_S (f_2 \circ g_2))\sum_i(m_i \otimes_S n_i) \nonumber \\
&=& \sum_i ((f_1 \circ g_1)(m_i) \otimes_S (f_2 \circ g_2)(n_i))\nonumber \\
&=& (f_1 \otimes_S f_2) \sum_i (g_1(m_i) \otimes_S g_2(n_i)) \nonumber \\
&=& (f_1 \otimes _S f_2) \circ (g_1 \otimes_S g_2) \sum_i (m_i \otimes_S n_i) \nonumber \\
&=& \otimes_S (f) \circ \otimes_S (g) \sum_i (m_i \otimes_S n_i); \nonumber \\
&&    
\end{eqnarray}
where the last equation holds whenever $f\circ g$ is defined in $(R,S) \times (S,T)$.  
Hence $\otimes_S$ is a functor.
However, we will see later that $\otimes_S$ is not associative. 
\label{er=tens}
\end{subsection}

\begin{subsection}{Traditional Morita theory}
Following~\cite{lam}, we will discuss the traditional Morita theory. Our approach emphasizes the algebraic side of the theory, starting with modules and generators. 
Later on, categories and functors will appear. However, see~\cite{Cohn} for a review of Morita theory which stresses the functoriality.
All theory below concerns right modules. Of course, an equivalent theory for left modules exists. 
First, we need some general notions and definitions.

\begin{definitie}
Let $R$ be a ring. Then $\mathfrak{M}_{R}$ denotes 
the category of right $R$-modules, the arrows of $\mathfrak{M}_{R}$ being given 
by $R$-module maps.
\end{definitie}

\begin{definitie} Let $R,S$ be two rings. $R$ and $S$ are called \textdef{Morita equivalent}, denoted by $R \meq S$, if there exists a 
categorical equivalence between $\mathfrak{M}_R$ and $\mathfrak{M}_S$, i.e. a functor $F : \mathfrak{M}_R \rightarrow \mathfrak{M}_S$ and a functor 
$ G : \mathfrak{M}_S \rightarrow \mathfrak{M}_R$ such that $(F \circ G) \simeq id_{\mathfrak{M}_S}$ and $(G \circ F) \simeq id_{\mathfrak{M}_R}$.
\label{defmorring}
\end{definitie}

\begin{definitie} A right $R$-module $P$ is a \textdef{generator} for $\mathfrak{M}_R$ if ${\mathrm{Hom}}_R(P,-)$ is a faithful functor from $\mathfrak{M}_R$ to 
the category of abelian groups. A finitely generated projective generator is called a \textdef{progenerator}. An \textdef{$(R,S)$-progenerator} $P$ is a faithfully balanced $R$-$S$ bimodule (i.e. a bimodule for which $R \cong {\mathrm{End}}(P_S)$ and $S \cong {\mathrm{End}}(_RP)$) that is a progenerator for $\mathfrak{M}_S$.
\end{definitie}

Recall that a module $P$ is finitely generated when for all families of submodules $\{ N_i\}_{i \in I}$ with $\sum_{i \in I} N_i = P$ there exists a 
finite subset $J \subseteq I$ such that $\sum_{i \in J} N_i = P$. Projectivity of $P$ implies that for any epimorphism (i.e. surjective homomorphism) 
of right $R$-modules $g : M \rightarrow N$  
and any $R$-homomorphism $h : P \rightarrow N$ there exists an $R$-homomorphism $f$ such that $h = g \circ f$.
\begin{equation}
\xymatrix{
& P\ar[d]^{h} \ar@{.>}[ld]_{f} & \\
M \ar[r]^{g} & N \ar[r] & 0.     \\
}
\end{equation}
The most trivial example of a generator for $\mathfrak{M}_R$ is the right regular module $R_R$. The functor ${\mathrm{Hom}}_R(R_R,-)$ is the forgetful functor from 
 $\mathfrak{M}_R$ to the category of abelian groups, which is faithful. Since $R_R$ is free, it is a projective module as well, because the universal 
property of a free module guarantees that a lifting as above exists. Hence $R_R$ is a progenerator.  
Note that being a (pro)generator can be expressed in terms of arrows and objects instead of elements, so being a (pro)generator is a categorical property. These facts will be used later in this section.

The following lemmas lead to a construction called ``Morita-context'', which is the basis for the Morita theorems.
\begin{lemma}
Let $P \in \mathfrak{M}_R$, and write  $P^* = {\mathrm{Hom}}_R(P,R)$ Consider the ring ${\mathrm{End}}_R(P)$. 
The following holds:
\begin{enumerate}
\item{
$P$ is an ${\mathrm{End}}_R(P)$-$R$ bimodule.}
\item{
$P^*$ is an $R$-${\mathrm{End}}_R(P)$ bimodule.}
\item{
There exists an $R$-$R$ homomorphism $\alpha : P^* \otimes_{{\mathrm{End}}_R(P)} P \rightarrow R$.\label{alphadef}}
\item{
There exists an  $S$-$S$ homomorphism $\beta : P \otimes_R P^* \rightarrow {\mathrm{End}}_R(P)$.\label{betadef}}
\end{enumerate}
\begin{bewijs}
Trivially,
${\mathrm{End}}_R(P)$ acts on the left of $P$ by applying the elements of ${\mathrm{End}}_R(P)$ to $P$. Recall that $P^*$ contains arrows $ P \leftarrow R$.
The right ${\mathrm{End}}_R(P)$-action on $P^*$ is defined by $(fg)(p) := f(gp)$, for $f \in P^*, g \in {\mathrm{End}}_R(P), p \in P$,
 which makes sense because of the left ${\mathrm{End}}_R(P)$-action on $P$. 
Viewing $R$ as $_RR_R$ it follows that $P^*$ has a left action of $R$ as well. Define $(rf)(p) := r(fp), r \in R, f \in P^*, p \in P$. 
This proves the first two parts of the lemma. 

Proof of \ref{alphadef}: Note that the fact that $P$ is an ${\mathrm{End}}_R(P)$-$R$ bimodule and  $P^*$ is an $R$-${\mathrm{End}}_R(P)$ bimodule implies that the 
tensor product $P^* \otimes_{{\mathrm{End}}_R(P)} P$ can be formed, and is an $R$-$R$ bimodule, see Remark \ref{simpeltens}. Define a mapping $\alpha$ as follows: 
\begin{eqnarray}
\alpha : P^* \otimes_{{\mathrm{End}}_R(P)} P &\rightarrow& R, \nonumber \\
   (f \otimes_{{\mathrm{End}}_R(P)} p) &\mapsto&  f(p),
\end{eqnarray}
for $f \in P^*, p \in P$. 
This $\alpha$ is well-defined: An easy computation shows that $\alpha$ vanishes 
on all elements $\sum_i(f_i \otimes_{{\mathrm{End}}_R(P)} p_i) = 0$. For example, 
for $f \in P^*, g \in {\mathrm{End}}_R(P), p \in P$ we have $(fg,p) = (f,gp)$. Now 
\begin{equation}
(fg,p) \mapsto fg(p) = f(gp),
\end{equation}
and 
\begin{equation}
(q,gp) \mapsto f(gp).
\end{equation}
The fact that $\alpha$ is a homomorphism follows from the fact that $P^*$ consists of homomorphisms. 
The $R$-$R$ action is preserved, since $P^*$ has a left $R$-action and $P$ has a right $R$-action.

Proof of \ref{betadef}: As above, the tensor product $P\otimes_R P^*$ can be formed. Define 
\begin{eqnarray}
\beta \ :\  P \otimes_R P^* &\rightarrow& {\mathrm{End}}_R(P),
\end{eqnarray}
by
\begin{equation}
\beta(p \otimes_R f)p^{\prime} = p(f(p^{\prime})),  
\label{pqp}
\end{equation}
$f \in P^*, p, p^{\prime} \in P$.
Note that in equation \eqref{pqp}, $f(p^{\prime})  \in R$, so $p(f(p^{\prime}))$ is defined by the right $R$-action on $P$.
This way, $\beta(p\otimes_R f)$ is an $R$-endomorphism of $P$: 

\begin{tabbing}
$\beta(p\otimes_Rf)(p_1 + p_2)$ \= = $p(f(p_1 + p_2))$ \hspace{80pt}\= (by definition of\ $\beta$) \\
\> = $p(f(p_1) + f(p_2))$ \> ($f$ is a homomorphism)  \\
\> = $p(f(p_1)) + p(f(p_2))$ \> ($P$ is a right $R$-module)  \\
\> = $\beta(p\otimes_R f)(p_1) + \beta(p\otimes_R f)(p_2)$ \> (by definition of\ $\beta$). \\ 
\text{Further,} \\
$\beta(p\otimes_R f)(p^{\prime} r)$ \> = $p(f(p^{\prime} r))$  \>  (by definition of\ $\beta$)  \\
\> = $p(f(p^{\prime})r)$ \> ($f$ is $R$-linear)  \\
\> = $p(f(p^{\prime}))r$ \> ($P$ is a right $R$ module)  \\
\> = $\beta((p\otimes_R f)(p^{\prime}))r$ \> (by definition of\ $\beta$).\\
\end{tabbing}
Similarly to the proof of \ref{alphadef}, one can show that $\beta$ is an 
${\mathrm{End}}_R(P)$-${\mathrm{End}}_R(P)$ homomorphism.
\end{bewijs} 
\label{beta=hom}
\end{lemma}

\begin{definitie} In the above notation, the $6$-tuple $(R,P,P^*,{\mathrm{End}}_R(P);\alpha,\beta)$ is called the \textdef{Morita-context} associated 
with $P_R$. The ring $R$ is called the \textdef{ground ring} for the Morita-context. 
More generally, for any ring $T$ and any right $T$-module $L$, the Morita-context associated with $L_T$ is 
\begin{equation}
(T,L,{\mathrm{Hom}}_T(L,T),{\mathrm{End}}_T(L);\alpha,\beta),
\end{equation}
where in this case $\alpha$ and $\beta$ are defined following the proof of \eqref{alphadef} and \eqref{betadef}.
\end{definitie}

Let $P_R$ be a right $R$-module, and fix the Morita-context \\ 
$(R,P,P^*,{\mathrm{End}}_R(P);\alpha,\beta)$. The following lemmas will show the 
connection between the notion of (pro)generator and the Morita-context. Despite the fact that we will use the results extensively, 
we will omit the proofs, since the lengthy computations do not provide much insight.
The proofs can be found in~\cite{lam}, Props. 18.17 and 18.19.

\begin{lemma}
\begin{enumerate} 
\item{
$P_R$ is a generator for $\mathfrak{M}_R$ iff $\alpha$ is an epimorphism.}
\item{
If one (and thus both) of the above conditions hold, then
\begin{enumerate}
\item{
$ \alpha : P^* \otimes_{{\mathrm{End}}_R(P)} P \rightarrow R$ is an $R$-$R$ isomorphism.
\label{linksend}
}
\item{
$ P^* \cong {\mathrm{Hom}}_{{\mathrm{End}}_R(P)}(P,{\mathrm{End}}_R(P))$ as $R$-${\mathrm{End}}_R(P)$ bimodules.
\label{rechtsend}
}
\item{
$ P \cong {\mathrm{Hom}}_{{\mathrm{End}}_R(P)}(P^*,{\mathrm{End}}_R(P))$ as ${\mathrm{End}}_R(P)$-$R$ bimodules.
}
\item{
$R \cong {\mathrm{End}}(_{{\mathrm{End}}_R(P)}P) \cong {\mathrm{End}}(P^*_{{\mathrm{End}}_R(P)})$ as rings.
}
\end{enumerate}
}
\end{enumerate}
\label{genlem}
\end{lemma}
One should note that the elements of the homomorphism space \\ 
\hbox{ ${\mathrm{Hom}}_{{\mathrm{End}}_R(P)}(P,{\mathrm{End}}_R(P))$} in \eqref{linksend} intertwine the \textdef{left} ${\mathrm{End}}_R(P)$-action, 
and that the elements of the homomorphism space \hbox{ ${\mathrm{Hom}}_{{\mathrm{End}}_R(P)}(P^*,{\mathrm{End}}_R(P))$} in \eqref{rechtsend} intertwine the \textdef{right} ${\mathrm{End}}_R(P)$-action.

\begin{lemma}
\begin{enumerate} 
\item{
$P_R$ is a finitely generated projective module iff $\beta$ is an epimorphism.}
\item{
If one (and thus both) of the above conditions hold, then
\begin{enumerate}
\item{
$ \beta : P \otimes_R P^* \rightarrow {\mathrm{End}}_R(P)$ is an ${\mathrm{End}}_R(P)$-${\mathrm{End}}_R(P)$ 
isomorphism.
}
\item{
$ P^* \cong {\mathrm{Hom}}_R(P_R,R_R)$ as $R$-${\mathrm{End}}_R(P)$ bimodules.
\label{!!}}
\item{
$ P \cong {\mathrm{Hom}}_R(_RP^*,_R\!R)$ as ${\mathrm{End}}_R(P)$-$R$ bimodules.
}
\item{
${\mathrm{End}}(P_R) \cong {\mathrm{End}}(_RP^*)$ as rings.
}
\end{enumerate}
}
\end{enumerate}
\label{proglem}
\end{lemma}
Note that the isomorphism in part \eqref{!!} is the identity by definition, it 
is included for symmetry reasons. 
\begin{stelling}[Morita I]
Let $P_R$ be a progenerator for $\mathfrak{M}_R$ and let\\
$(R,P,P^*,{\mathrm{End}}_R(P);\alpha,\beta)$ be the Morita-context associated with $P_R$. 
Then the functors 
\begin{equation}
-\otimes_R P^* : \mathfrak{M}_R \rightarrow 
\mathfrak{M}_{{\mathrm{End}}_R(P)}
\end{equation}
 and 
\begin{equation}
-\otimes_{{\mathrm{End}}_R(P)} P : \mathfrak{M}_{{\mathrm{End}}_R(P)} 
\rightarrow \mathfrak{M}_R
\end{equation}
 are mutually inverse (up to natural isomorphism) and hence they are 
category equivalences.
\begin{bewijs}
Let $U_R \in \mathfrak{M}_R$. Then 
\begin{eqnarray}
(U \otimes_R P^*) \otimes_{{\mathrm{End}}_R(P)} P &\simeq& U \otimes_R ( P^* \otimes_{{\mathrm{End}}_R(P)} P) \nonumber \\ &\simeq& 
U \otimes_R R \hspace{10pt}\simeq \hspace{10pt}U,
\label{objiso}
\end{eqnarray}
where the first isomorphism follows from the fact that the tensor product is 
associative up to isomorphism (see Proposition \ref{rings=bic}). The 
 second isomorphism follows from the fact that $P_R$ 
is a generator, so that $\alpha$ is an isomorphism (using Lemma \ref{genlem}). 
The other
isomorphism follows from the construction of the tensor product. 
Using \eqref{objiso} we have that on objects
\begin{equation}
(- \otimes_{{\mathrm{End}}_R(P)}P) \circ (-\otimes_R P^*) \simeq id_{\mathfrak{M}_R}.
\end{equation}
Now for $f : U_R \rightarrow V_R$ in $(\mathfrak{M}_R)_1$, the diagram
\begin{equation}
\xymatrix{
(U \otimes_R P^*) \otimes_{{\mathrm{End}}_R(P)} P \ar[rr] \ar[d]_{(f \otimes_R \ id_{P^*}) 
\otimes_{{\mathrm{End}}_R(P)} id_P} && U\ar[d]^f \\
(V \otimes_R P^*) \otimes_{{\mathrm{End}}_R(P)} P \ar[rr] && V,
}
\end{equation}

commutes by definition of the tensor product on arrows.

Let $W_{{\mathrm{End}}_R(P)} \in \mathfrak{M}_{{\mathrm{End}}_R(P)}$. Then, similarly, 
\begin{eqnarray}
(W \otimes_{{\mathrm{End}}_R(P)} P) \otimes_R P^* &\simeq& W \otimes_{{\mathrm{End}}_R(P)} (P \otimes_R P^*) \nonumber \\ &\simeq& 
W \otimes_{{\mathrm{End}}_R(P)} {\mathrm{End}}_R(P) \hspace{8pt}\simeq \hspace{8pt}W,
\end{eqnarray}
and a similar commutative diagram as above can be constructed.
\end{bewijs}
\end{stelling}

\begin{lemma}
Let $P_R$ be a progenerator for $\mathfrak{M}_R$, with associated Morita-context $(R,P,P^*,{\mathrm{End}}_R(P);\alpha,\beta)$. One has 
the following natural functor isomorphisms:
\begin{enumerate}
\item{
$- \otimes_R P^* \simeq {\mathrm{Hom}}_R(P_R,-),$
}
\item{
$- \otimes_{{\mathrm{End}}_R(P)} P \simeq {\mathrm{Hom}}_{{\mathrm{End}}_R(P)}(P^*_{{\mathrm{End}}_R(P)},-),$}
\item{
$P \otimes_R -\simeq {\mathrm{Hom}}_R(_R\!P^*,-),$}
\item{
$
P^* \otimes_{{\mathrm{End}}_R(P)} - \simeq {\mathrm{Hom}}_{{\mathrm{End}}_R(P)}(_{{\mathrm{End}}_R(P)}P,-).$}
\end{enumerate}
\begin{bewijs}
We construct only the first functor isomorphism. The other parts of the lemma can be proved likewise.
Take $M_R \in (\mathfrak{M}_R)_0$. Define 
\begin{eqnarray}
\beta_M : M \otimes_R P^* &\longrightarrow& {\mathrm{Hom}}_R(P_R,M_R), \nonumber \\
\text{by}\hspace{150pt} & \nonumber \\
\beta_M(m \otimes_R f)(p) &=&  m(fp). 
\end{eqnarray}
Defined like this, $\beta_M$ is an isomorphism of right ${\mathrm{End}}_R(P)$-modules. 
The fact that $\beta_M$ is a homomorphism is proven analogously to Lemma \ref{beta=hom}. 
Injectivity is trivial, and surjectivity follows from that fact that $P_R$ is a projective module. $\beta_M$ intertwines the ${\mathrm{End}}_R(P)$-action by definition of the 
right ${\mathrm{End}}_R(P)$-action on $P^*$. 
An isomorphism like this 
can be formed for all right $R$-modules in $\mathfrak{M}_R$. Therefore, we can identify the functor $-\otimes_R P^*$ with 
the functor ${\mathrm{Hom}}_R(P_R,-)$. For $g : M \rightarrow N \in (\mathfrak{M}_R)_1$, 
we will show that the diagram 
\begin{equation} 
\xymatrix{
M \otimes_R P^* \ar[rr]^{\beta_M} \ar[d] & & {\mathrm{Hom}}_R(P_R,M_R)\ar[d] \\
N \otimes_R P^* \ar[rr]^{\beta_N} & & {\mathrm{Hom}}_R(P_R,N_R) \\
}
\end{equation}
is commutative. In the upper half of the diagram $$ (m,f) \mapsto m(f(-)) 
\mapsto g(m(f(-)))$$ holds. In the lower half of the diagram we have $$(m,f) 
\mapsto (g(m),f) \mapsto g(m)(f(-)) = g(m(f(-))),$$ since $g$ intertwines 
$f(-) \in R$.   
\end{bewijs}
\label{homiso}
\end{lemma}
Note the similarity of the lemma above to some of the statements of Lemmas \ref{genlem} and \ref{proglem}. For example, taking $M = P$ in the proof above, 
we obtain $\beta_M = \beta : P\otimes_R P^* \rightarrow \mathrm{Hom}_R(P_R,P_R) = \mathrm{End}_R(P)$, so we may consider the functor isomorphism above as a generalization 
of $\beta$.
\begin{stelling}[Morita II]
Let $R,S$ be two rings, and 
$$ F : \mathfrak{M}_R \rightarrow \mathfrak{M}_S, 
\hspace{30pt} G : \mathfrak{M}_S \rightarrow \mathfrak{M}_R $$
be mutually inverse category equivalences. 
Consider $ F(R_R)$ and $G(S_S)$. We have functor isomorphisms, 
i.e. natural equivalences 
\begin{eqnarray}
F &\simeq& - \otimes_R F(R_R)  \hspace{30pt} \text{and} \label{1stdef} \\ 
G &\simeq& -  \otimes_S G(S_S). \label{2nddef}
\end{eqnarray}
\begin{bewijs}
First, we show that $F(R_R)$ and $G(S_S)$ have bimodule structures:\\
 \hbox{$F(R_R) = _R F(R_R)_S $} and \hbox{$G(S_S) = _S G(S_S)_R$.} 
By definition of $F$ and $G$, $F(R_R)$ has 
a right $S$-action and $G(S_S)$ has a right $R$-action. Further, $F(R_R)$ inherits a right $R$-action from $_RR$ through $F_1$, since 
$R \subset (R,R) \subset (\mathfrak{M}_R)_1$. We can define $rk := F_1(r)k$, $r \in R, k \in F(R_R)$.
Hence $F(R_R)$ can be seen as an $R$-$S$ bimodule. A similar 
argument holds for $G(S_S) = _S G(S_S)_R$. This shows 
that the tensor products in 
\eqref{1stdef} and \eqref{2nddef} can indeed be formed.

Further, since being a progenerator is a categorical property, $F(R_R)_S$ 
is a progenerator for $\mathfrak{M}_S$ because $R_R$ is for 
$\mathfrak{M}_R$. Now we compute 
 $F(R_R)^*$ as follows:.
\begin{eqnarray}
F(R_R)^* &=& {\mathrm{Hom}}_S(S,F(R_R)) \simeq {\mathrm{Hom}}_R(G(S),G(F(R_R))) 
\nonumber \\ &\simeq& {\mathrm{Hom}}_R(G(S_S)_R,R_R) \simeq G(S_S).
\end{eqnarray} 
Therefore the Morita-context associated with the progenerator $F(R_R)_S$ is 
\\$(S,F(R_R)_S,G(S_S),R;\alpha,\beta)$, where $\alpha$, $\beta$ are the 
appropriate pairings. In particular, the Morita I can 
be applied.

To show the natural equivalence \eqref{2nddef}, let $N_S 
\in \mathfrak{M}_S$. Now, 
\begin{equation}
G(N_S) \simeq {\mathrm{Hom}}_S(G(N_S),R_R) \simeq {\mathrm{Hom}}_R(N_S,F(R_R)_S).
\end{equation}
Applying Lemma \ref{homiso} proves \eqref{2nddef}. The other part of the theorem,
 $$F \simeq {\mathrm{Hom}}_R(-, G(S_S)_R) \simeq - \otimes_R F(R_R), $$ 
follows by a similar argument.  
\end{bewijs}
\end{stelling}

Observe that Morita II is a special case of the Eilenberg-Watts theorem, 
 which states that (a certain class of) functors between $\mathfrak{M}_R$ and 
$\mathfrak{M}_S$ are natural equivalent to taking tensor products with an 
appropriate bimodule, see~\cite{bass} or~\cite{watts}. Originally, the Eilenberg-Watts theorem holds for 
algebras. Viewing rings as algebras over $\mathbb{Z}$, it is applicable to rings as well.
\begin{opmerking}

Following the same line of argument as above for $G(S_S)_R$, one may define a
 Morita-context 
$(R,G(S_S)_R,F(R_R),S;\alpha,\beta)$ associated with the progenerator 
$G(S_S)_R$ of $\mathfrak{M}_R$. The situation in 
Morita II is symmetric in $R$ and $S$. Given a pair of 
mutually inverse equivalences between $\mathfrak{M}_R$ and $\mathfrak{M}_S$ 
one can construct a Morita-context 
with ground ring $R$ or $S$. Note the difference with the first 
Morita theorem where the Morita-context is fixed: there we obtain a pair category equivalences from the rings $R$ and ${\mathrm{End}}_R(P)$. 
\end{opmerking}

The following corollary is the main result of Morita theory. It immediately follows from Morita I and II. 
\begin{gevolg}
Two rings are 
Morita equivalent iff one is isomorphic to the endomorphism ring of some progenerator of the other 
(and vice versa).  
\end{gevolg}
\newpage
\begin{stelling}[Morita III]
Let $R$,$S$,$T$ be rings. Then a bijective correspondence exists between the 
isomorphism classes of the category equivalences $\mathfrak{M}_S 
\rightarrow \mathfrak{M}_R$ and the isomorphism classes of $(S,R)$-progenerators.

Composition of category equivalences $\mathfrak{M}_T \rightarrow \mathfrak{M}_S 
\rightarrow \mathfrak{M}_R$ corresponds to tensor products of $(T,S)$ and $(S,R)$
-progenerators.
\begin{bewijs}
Recall that an $(S,R)$-progenerator $_SP_R$ is an $S$-$R$ bimodule $P$ which 
is a progenerator for $\mathfrak{M}_R$ and $R \cong {\mathrm{End}}(_SP), S \cong {\mathrm{End}}(P_R)$. 
A $(S,R)$-progenerator $_SP_R$ leads to a category equivalence 
$-\otimes_S P : \mathfrak{M}_S \rightarrow \mathfrak{M}_R$. We have already seen that $-\otimes_S P$ 
is a functor; Morita I states that it is a category 
equivalence as well, using $S \cong {\mathrm{End}}(P_R)$. 
For an $(S,R)$-progenerator $_SQ_R$ that is isomorphic 
to $_SP_R$ via an isomorphism $\gamma$, we construct a category equivalence 
$-\otimes_S Q$. This category equivalence is in the same isomorphism class 
as $- \otimes_S P$ since 
\begin{equation}
id \otimes_S \gamma : K \otimes_S P \rightarrow K \otimes_S Q,
\end{equation}
is a natural isomorphism.

Conversely, suppose we have a category equivalence $F : 
\mathfrak{M}_S \rightarrow \mathfrak{M}_R$. Then $F(S_S)$ is a 
progenerator for $\mathfrak{M}_R$, since $S_S$ is a progenerator 
for $\mathfrak{M}_S$. As we have seen in the proof of Morita II, $F(S_S)$ is 
an $S$-$R$ bimodule.
Let $G: \mathfrak{M}_S \rightarrow \mathfrak{M}_R$ be a 
category equivalence in the same isomorphism class as $F$. Then $G(S_S)$ is a 
progenerator for $\mathfrak{M}_S$ and an $S$-$R$ bimodule as well. Further, 
\begin{equation}
G(S_S) \simeq F(S_S),
\end{equation} 
since $G \simeq F$.
This proves the bijective correspondence in the theorem. 

Let $P$ be an $(S,R)$-progenerator and $Q$ be a $(T,S)$-progenerator. Then 
$P$, $Q$ give rise to category equivalences
\begin{equation}
- \otimes_S P : \mathfrak{M}_S \rightarrow \mathfrak{M}_R,
\end{equation}
\begin{equation}
- \otimes_T Q : \mathfrak{M}_T  \rightarrow \mathfrak{M}_S.
\end{equation}
Composition of these category equivalences gives the category equivalence
\begin{equation}
(- \otimes_T Q) \otimes_S P : \mathfrak{M}_T \rightarrow \mathfrak{M}_R,
\end{equation}
which is isomorphic (via natural equivalence) to 
\begin{equation}
- \otimes_T (Q \otimes_S P) : \mathfrak{M}_T \rightarrow \mathfrak{M}_R.
\end{equation}
The last category equivalence is just the equivalence we obtain when we first 
apply the tensor product to the progenerators $Q$ and $P$, and then generate 
a category equivalence from the tensor product $Q \otimes_S P$.
\end{bewijs}
\end{stelling}
\end{subsection}

\begin{subsection}{Morita theory with use of bicategories}

To state the Morita theory in terms of bicategories, we first need to show that 
rings, bimodules and bilinear maps indeed form a bicategory. This fact (without proof) was already stated in~\cite{ben} and~\cite{maclane}. Recall Section \ref{bic}.

\begin{prop}
For any two rings $R$,$S$, let $(R,S)$ be the category of $R$-$S$ bimodules as objects, and $R$-$S$ linear maps as arrows.
Then the collection of all rings as objects and bimodules as arrows forms a bicategory $\mathsf{[Rings]}$, in which the composition functor 
$(R,S) \times (S,T) \rightarrow (R,T)$ is given by the tensor product $\otimes_S$, and the unit arrow in $(R,R)$ is given by 
$I_{R} = R \rightarrow R \leftarrow R$.
\begin{bewijs}
One by one, the properties of a bicategory as stated in Section \ref{bic} will be checked.
\begin{itemize}
\item{
The object space $\mathsf{[Rings]}_0$ consists of all rings.}
\item{
For any two rings, there is a category $(R,S)$. $(R,S)$ has $R$-$S$ bimodules as objects and $R$-$S$ linear maps as arrows.
It is easy to see that $(R,S)$ is indeed a category, we will not prove this here.}
\item{
For each triple $(R,S,T)$ of rings there is a composition functor given by the tensor product $\otimes_S$. 
\begin{eqnarray}
 \otimes_S : (R,S) \times (S,T) &\longrightarrow& (R,T), \nonumber \\
 _RM_S \times _SN_T &\longmapsto& _R(M \otimes_S N)_T,
\end{eqnarray}
for $ _RM_S \in (R,S), _SN_T \in (S,T)$.
It has been shown in Section \ref{er=tens} that $\otimes_S$ is indeed a functor, and that it is properly defined.
}
\item{For each object $R$ of $\mathsf{[Rings]}_0$, the identity arrow $I_R$ of $(R,R)$ is given by $R\rightarrow R \leftarrow R$.
}
\item{For each quadruple $(R,S,T,U)$ of rings, 
we need a natural isomorphism $\beta$ between the functors $
((- \otimes_S - ) \otimes_T- ) $ and $ (- \otimes_S (- \otimes_T -))$, each from the category  
$(R,S) \times (S,T) \times (T,U)$ to the category $(R,U)$. To each object $(M,N,P)$ in $(R,S) \times (S,T) \times (T,U)$, let $\beta$ assign an arrow 
$$((M \otimes_S N) \otimes_T P) \rightarrow (M \otimes_S ( N \otimes_T P)),$$ in $(R,U)$, where $\beta$ is defined as follows on the elements: $$(m\otimes_S n) \otimes_T p 
\mapsto m \otimes_S (n \otimes_T p).$$
Now, given an arrow $$(\phi,\psi,\chi) : (M,N,P) \rightarrow (M^{\prime},N^{\prime},P^{\prime})$$
in $(R,S) \times (S,T) \times (T,U)$, the diagram
\begin{equation}
\xymatrix{
 (M \otimes_S N) \otimes_T P \ar[d]^{(\phi \otimes_S \psi) \otimes_T \chi} 
\ar[rrr]^{\beta(M,N,P)} &&& M \otimes_S (N \otimes_T P ) \ar[d]^{\phi \otimes_S (\psi \otimes_T \chi)} \\
\bigl(M^{\prime} \otimes_S N^{\prime} \bigr) \otimes_T P^{\prime} \ar[rrr]^{\beta(M^{\prime},N^{\prime},P^{\prime})} &&&  M^{\prime} \otimes_S \bigl(N^{\prime} \otimes_T P^{\prime} \bigr)\\
}
\end{equation}
commutes by definition of the tensor functor on arrows.
}
\item{For each pair of rings $(R,S)$, we need the left identity $L_{(R,S)}$, which is a natural isomorphism between the functor $\otimes_R \circ \left[I_R \times Id_{(R,S)}\right]$ and the canonical 
functor from $1 \times (R,S)$ to $(R,S)$. Let $L_{(R,S)}$ assign an arrow $ R \otimes_{R} M_S  \rightarrow _{R}\!M_S$ in $(R,S)$ to each object in $ 1 \times (R,S)$. Then
$L_{(R,S)}( 1 \times _{R}\!M_S)$ is given by
\begin{eqnarray}
L_{(R,S)}( 1 \times _{R}\!M_S) : R \otimes_R M_S  &\rightarrow& _{R}M_S \nonumber \\
(r \otimes_R m) &\mapsto& rm. 
\end{eqnarray}
This is an invertible arrow, by sending $m \mapsto (1,m)$, since $(r,m) = (1,rm)$ in $R\otimes_R M$.
Now, given an arrow $(1,\phi) : 1 \times _RM_S \rightarrow 1 \times _RN_S$ in $1 \times (R,S)$, the diagram 
\begin{equation}  
\xymatrix{   
R \otimes_R M \ar[d]^{1\otimes_R \phi} \ar[rrr]^{L_{(R,S)}(1 \times _RM_S)} &&& _RM_S\ar[d]^{\phi} \\
R \otimes_R N \ar[rrr]^{L_{(R,S)}(1 \times _RN_S)} &&& _RN_S, \\}
\end{equation}
is commutative. Take $(r,m) \in R \otimes_R M$. Then \hbox{$ (r,m) \mapsto rm \mapsto \phi(r,m)$} in the upper half of the diagram. In the lower half of the 
diagram
\hbox{$(r,m) \mapsto (r,\phi(m)) \mapsto r\phi(m) = \phi(rm)$} where we used that $\phi$ is $R$-linear.

The right identity $R_{(R,S)}$ is defined likewise.
}
\item{
We need to prove associativity coherence. Let $R,S,T,U,V$ be rings, and $_RK_S,\ _SL_T,\ _TM_U,\ _UN_V$ be associated bimodules. The diagram
\begin{equation}
\xymatrix{
((K \otimes_S L) \otimes_T M) \otimes_U N \ar[rrr]^{\beta(K,L,M) \otimes_U Id}\ar[d]_{\beta(K\otimes_SL,M,N)} &&&
(K \otimes_S (L \otimes_T M)) \otimes_U N \ar[d]^{\beta(K,L\otimes_TM,N)} \\
(K \otimes_S L) \otimes_T (M \otimes_U N) \ar[dd]_{\beta(K,L,M\otimes_UN)} &&& K \otimes_S ((L \otimes_T M ) \otimes_U N) \ar[ddlll]^{\ \ \ \ \ Id \otimes_S \beta(L,M,N)} \\ &&&& \\
 K \otimes_S (L \otimes_T (M \otimes_U N)), &&\\
}
\end{equation}
commutes, since $\beta$ is defined elementwise.
}
\item{
We need to prove identity coherence. Let $R,S,T$ be rings with associated bimodules $_RP_S, _SQ_T$. We show that the following diagram commutes: 
\begin{equation}
\xymatrix{
(P \otimes_S S) \otimes_S Q \ar[rr]^{\beta(P,S,Q)}\ar[dr]_{R(P,S) \otimes_S Id\ \ \ } && P \otimes_S ( S \otimes_S Q)\ar[dl]^{\ \ \ \ \ Id \otimes_S L(S,Q)} \\
& P\otimes_S Q. & \\
}
\label{diag}
\end{equation}
Let $(p \otimes_S s) \otimes_S q$ be an element of $ (P \otimes_S S) \otimes_S Q$. In the upper half of the diagram we have
$$ (p \otimes_S s) \otimes_S q \mapsto p \otimes_S (s \otimes_S q) \mapsto p \otimes_S sq. $$
In the lower half of the diagram, the equality 
$$ (p \otimes_S s) \otimes_S q \mapsto ps \otimes_S q = p \otimes_S sq, $$
holds. Therefore, Diagram \eqref{diag} commutes. 
}
\end{itemize}
\end{bewijs}
\label{rings=bic}
\end{prop}

Recall Definition \ref{defmorring}; equivalence of the categories of right modules of two rings is called Morita equivalence. The next result formulates Morita 
equivalence in terms of bicategories. 
  
\begin{stelling}
Two rings $R$ and $S$ are isomorphic objects in the bicategory $\mathsf{[Rings]}$ iff their
 categories of right modules are equivalent. In formula:
\begin{equation}
R \biso S \Longleftrightarrow \mathfrak{M}_{R} \simeq \mathfrak{M}_{S}.
\end{equation}
\begin{bewijs}
``$\Rightarrow:$'' Suppose $R \biso S$. Being isomorphic in $\mathsf{[Rings]}$ means that there exists a bimodule 
$_{R}N_{S}$ in $(R,S)$ which is invertible up to isomorphism. Therefore there exists a bimodule 
$_{S}N^{-1}_{R}$ in $(S,R)$ such that
\begin{eqnarray} R \rightarrow N\otimes_S N^{-1} \leftarrow R &\cong& R \rightarrow R \leftarrow R \text{\ \ in\ \ } (R,R); \label{invis1}\\
 S \rightarrow N^{-1}\otimes_{R} N \leftarrow S &\cong& S \rightarrow S \leftarrow S\text{\ \ in\ \ } (S,S). \label{invis2}
\end{eqnarray}
Now, a functor $F : \mathfrak{M}_{R} \longrightarrow \mathfrak{M}_S$ will be constructed. On objects, define 
\begin{equation}
F_0(K):= K \otimes_{R} N, \text{\ for\ } K \in (\mathfrak{M}_{R})_0.
\end{equation}
On arrows, one defines 
\begin{equation}
F_1(f) := f \otimes_{R} id_N, \text{\ for\ } f \in (\mathfrak{M}_{R})_1.
\end{equation}
This way, $F$ is a functor:
To each object $K$ in $\mathfrak{M}_{R}$, $F$ assigns an object $K \otimes_{R} N$ in $\mathfrak{M}_S$. $F_0(K)$ has a right action of $S$ 
 which is passed to the tensor product from the right $S$-action of $_{R}N_{S}$, see Section \ref{er=tens}. Further,
$F$ assigns to each arrow $h : K_1 \rightarrow K_2$ in $\mathfrak{M}_{R}$ an arrow $F_1(h) : K_1 \otimes_{R} N \rightarrow 
K_2 \otimes_{R} N$ in $\mathfrak{M}_S$.  Moreover,
\begin{eqnarray}
F_1(id_{K})\ (k \otimes_{R} n) &=& (id_{K} \otimes_{R} id_N)\  (k \otimes_{R} n)  \nonumber \\
&=& (id_K (k) \otimes_{R} id_N (n)) = k \otimes_{R} n\nonumber \\ &=& 
 id_{(K \otimes_{R} N)}(k \otimes_{R} n)  \nonumber \\ &=& id_{F_1(id_K)}, 
\end{eqnarray}
and
\begin{eqnarray}
F_1(h_1 \circ h_2)(k \otimes_{R} n) &=& (h_1 \circ h_2(k) \otimes_{R} n)  \nonumber \\ &=&
F_1(h_1)\ (h_2 (k) \otimes_{R} n)\nonumber\\ &=& F_1(h_1) \circ F_1 (h_2) \ (k \otimes_{R} n); 
\end{eqnarray}   
for $K \in (\mathfrak{M}_{R})_0$, $h_1$ and $h_2 \in (\mathfrak{M}_{R})_1$.
The last equation holds whenever $h_1 \circ h_2$ is defined in $(\mathfrak{M}_{R})_1$.

In the same way, one constructs a functor $G : \mathfrak{M}_S \longrightarrow \mathfrak{M}_{R}$ by putting 
$G_{0}(L) = L \otimes_S N^{-1}$ for $L \in (\mathfrak{M}_S)_0$ and $G_1(g) = g \otimes_S id_{N^{-1}}$ 
for $g \in (\mathfrak{M}_S)_1$. 

To prove equivalence, we need to show that natural equivalences \\$(F \circ G) 
\simeq id_{\mathfrak{M}_{S}}$ and $(G \circ F) \simeq id_{\mathfrak{M}_R}$ exist. Using \eqref{invis2} and the fact that the 
tensor product is associative up to isomorphism, we get 
\begin{eqnarray}
(F \circ G)_0(L)&=& F_0(L \otimes_S N^{-1}) \  = \  (L \otimes_S N^{-1}) \otimes_{R} N
\\ \nonumber &\simeq& L \otimes_S ( N^{-1} \otimes_{R} N) \ \simeq \ L \otimes_S S \hspace{8pt} \simeq \hspace{8pt}L, \label{bluh}
\end{eqnarray}
for $L \in (\mathfrak{M}_S)_0$. 

Let $g$ be an arrow $L_{S} \rightarrow M_{S}$ in $\mathfrak{M}_{S}$.
Clearly, the diagram,
\begin{equation}
\xymatrix{
F \circ G (L) \ar[rr] \ar[d]_{g \otimes_{S} (id_{N^{-1}} \otimes_R id_N)} && L \ar[d]^{g} \\
M \otimes_{S} (N^{-1} \otimes_R N)\ar[rr] && M, \\
}
\end{equation}
commutes, via the isomorphism \eqref{bluh}. A similar computation shows that \\
 $(G \circ F) \simeq id_{\mathfrak{M}_R}$. This proves the ``$\Rightarrow$'' part of the theorem.

To prove the ``$\Leftarrow$'' part of the theorem, we need to construct an invertible $R$-$S$ bimodule $N$. Given are two equivalent functors
 $F : \mathfrak{M}_{R} \longrightarrow \mathfrak{M}_S$ and $G : \mathfrak{M}_S \longrightarrow \mathfrak{M}_{R} $. Define $$ N = F_0(R_{R}). $$
As we have seen in the proof of Morita II, we have $N \in (R,S)$: By definition of $F_0$, $N$ has a right action of $S$. $N$ also inherits a left action of $R$ from the left action of $R$ on $R_{R}$, through 
$F_1$, since $R \subset (R,R) \subset (\mathfrak{M}_{R})_1$. Define 
$rn := F_1(r)n$. Thus $N \in (R,S)$. Similarly, one defines $N^{-1} = G_0(S_S)$. Following the same reasoning as above, $N^{-1} \in (S,R)$.
Since $F,G$ are mutually inverse category equivalences, we can apply Morita II, 
which immediately proves \eqref{invis1} and \eqref{invis2}. 
\begin{eqnarray}
N \otimes_S N^{-1} &=& F(R_R) \otimes_S G(S_S) \ \simeq \ G(F(R_R))\  \simeq \ R, \nonumber \\
N^{-1} \otimes_R N &=& G(S_S) \otimes_R F(R_R) \ \simeq \ F(G(S_S))\ \simeq \ S.
\end{eqnarray}
In other words, $N$ is an invertible $R$-$S$ bimodule, so $R \biso S$.
\end{bewijs}
\label{Bicat}
\end{stelling}
\label{Moritabic}
\end{subsection}

\begin{subsection}{Equivalence of theories}
This section will show equivalence between the ``traditional'' Morita theory 
and the Morita theory stated in terms of bicategories. In particular, we will prove 
equivalence between the Morita III theorem and Theorem \ref{Bicat}. The last 
theorem states that Morita equivalence of two rings in the traditional 
setting is equivalent to being isomorphic in the bicategory $\mathsf{[Rings]}$.
The proof of equivalence leans heavily on the following proposition.
\newpage
\begin{prop}
Let $_RM_S$ be an $(R,S)$-module. Then $_RM_S$ is an invertible module in $\mathsf{[Rings]}$ iff$^{\ }$\footnote{Due to the right-left symmetry of the theory, it is also 
possible to prove that $_RM_S$ is an invertible module iff 
\begin{enumerate}
\item{
$R_M$ is a progenerator for $_R{\mathfrak{M}}$.
}
\item{
$S \cong {\mathrm{End}}_{R}(M)$ as rings.
}
\end{enumerate}
Note that the second item follows immediately from Proposition \ref{crucial} and Lemma \ref{genlem}. To prove the first item, use 
$$ \mathfrak{M}_R \simeq \mathfrak{M}_S  \Leftrightarrow _R\!{\mathfrak{M}} \simeq _S\!{\mathfrak{M}}, $$
and a similar argument as in the proof of Proposition \ref{crucial}. See~\cite{lam} for further details.
}

\begin{enumerate}
\item{
$M_S$ is a progenerator for $\mathfrak{M}_S$.
}
\item{
$R \cong {\mathrm{End}}_{S^{op}}(M)$ as rings.
}
\end{enumerate}

\begin{bewijs}
$\Rightarrow$: Suppose the bimodule $R \rightarrow M \leftarrow S$ is invertible. Using Theorem \ref{Bicat}, a categorical 
equivalence $F : \mathfrak{M}_R \rightarrow \mathfrak{M}_S$ exists. Let $G : \mathfrak{M}_S \rightarrow \mathfrak{M}_R$ denote an inverse (up to natural equivalence).
As in the proof of Theorem \ref{Bicat}, an 
invertible module $R \rightarrow N \leftarrow S$ can be constructed by defining $N = F(R_R)$. $R_R$ is a progenerator for 
$\mathfrak{M}_R$, hence $N_S$ is a progenerator for $\mathfrak{M}_S$ since being a progenerator is a categorical property. 
The functor $F$ thus obtained acts on the objects of $\mathfrak{M}_R$ by $K_R \mapsto K \otimes_R M$, so we have 
\begin{equation} 
N_S = F(R_R) = R \otimes_R M \cong M_S.
\end{equation}    
This makes $M_S$ a progenerator for $\mathfrak{M}_S$.
Following Morita II, we obtain the Morita-context
$$ (F(R_R), S, G(S_S), R; \alpha, \beta),$$
for $F(R_R)$.
However, since $M_S \cong F(R_R)$ we have 
\begin{equation} 
R = {\mathrm{End}}_{S^{op}}(F(R_R)) \cong {\mathrm{End}}_{S^{op}}(M_S),
\end{equation}
which proves the ``$\Rightarrow$'' part of the proposition.

$\Leftarrow$: Let $M_S$ be a progenerator for $\mathfrak{M}_S$. Fix the Morita-context\\ $(S,M,{\mathrm{Hom}}_S(S,M),{\mathrm{End}}(M_S);\alpha,\beta)$. 
Using the fact that $R \cong {\mathrm{End}}(M_S)$ we obtain an $S$-$S$ isomorphism ${\mathrm{Hom}}_S(S,M) \otimes_R M \rightarrow S$ and an $R$-$R$ isomorphism  
$M \otimes_S {\mathrm{Hom}}_S(S,M) \rightarrow R$ from Lemmas \ref{genlem} and \ref{proglem}. So
\begin{eqnarray}
S \rightarrow {\mathrm{Hom}}_S(S,M) \otimes_R M \leftarrow S &\cong& S \rightarrow S \leftarrow S \\ 
R \rightarrow M \otimes_S {\mathrm{Hom}}_S(S,M) \leftarrow R &\cong& R \rightarrow R \leftarrow R 
\end{eqnarray}
This shows that $M_S$ is invertible, its inverse being equal to ${\mathrm{Hom}}_S(S,M)$. 
\end{bewijs}
\label{crucial}
\end{prop}
The remainder of this subsection is devoted to the proof that Theorem \ref{Bicat} can be derived from Morita III and vice versa. 
First we will show that Theorem \ref{Bicat} implies the Morita III. So assume Theorem \ref{Bicat}.
Let $R$,$S$ be two rings. 
Suppose we have an isomorphism class of category equivalences $\mathfrak{M}_R \rightarrow 
\mathfrak{M}_S$. Let $F$ be a representative in this class. According to Theorem \ref{Bicat}, 
 there exists an invertible $(R,S)$-module $_RM_S$, defined by $F(R_R)$. 
Applying Proposition \ref{crucial} it follows that $M_S$ is a progenerator for $\mathfrak{M}_S$. 

Let $\tilde{F} : \mathfrak{M}_R \rightarrow 
\mathfrak{M}_S$ be another representative in this isomorphism class. Hence we have a natural isomorphism $\sigma$ between $F$ and $\tilde{F}$. Once again,
$\tilde{F}(R_R)$ defines an invertible $(R,S)$-module $_RN_S.$ Since $F \simeq \tilde{F}$ through 
$\sigma$, we have that \begin{equation}
_RM_S = F(R_R) \simeq \tilde{F}(R_R) = _R\!N_S.
\end{equation} 
So elements of the isomorphism 
class of category equivalences lead to isomorphic progenerators. 

On the other hand, let $P_S,Q_S$ be two isomorphic $(R,S)$-progenerators for $\mathfrak{M}_S$. Proposition \ref{crucial} shows 
that $_RP_S$ and $_RQ_S$ are invertible modules. 
Following Theorem \ref{Bicat} we are able to construct two equivalence functors \\ $G,\tilde{G} : \mathfrak{M}_R 
\rightarrow \mathfrak{M}_S$, defined by taking the tensor product with $P$ and $Q$, respectively. 
Note that in the previous sentence, ``equivalence functors'' does \emph{not} mean that $G$ and $\tilde{G}$ are mutually 
inverse equivalences, for both $G$ and $\tilde{G}$ have the same domain! It only means that $G$ and $\tilde{G}$ each have an inverse. Recall that $G$ can be defined as follows: 
\begin{eqnarray} 
&&G : \mathfrak{M}_R \rightarrow \mathfrak{M}_S \nonumber \\
&&G_0(L) = L \otimes_R P \hspace{10pt} \text{for } L \in (\mathfrak{M}_R)_0 \nonumber \\
&&G_1(g) = g \otimes_R id \hspace{10pt} \text{for } g \in (\mathfrak{M}_R)_1. 
\end{eqnarray}
The functor $\tilde{G}$ is defined likewise.
It already has been shown that $G$ and $\tilde{G}$ are functors and have an inverse; see the proof of Theorem \ref{Bicat}.
We will construct a natural isomorphism $\tau$: To each object $L\in \mathfrak{M}_R$ let $\tau$ assign an invertible arrow $L\otimes_R P \mapsto L \otimes_R Q$. We can do this by defining $\tau$ as $id \otimes_R \rho$, where $\rho$ is given by the isomorphism between $_RP_S$ and $_RQ_S$.
 Let $\phi$ be an arrow $L \mapsto L^{\prime}$ in $\mathfrak{M}_R$. We obtain the 
following commutative diagram:
\begin{equation}
\xymatrix{
L \otimes_R P \ar[d]^{\phi \otimes_R id} \ar[rr]^{\tau(L)} && L\otimes_R Q \ar[d]^{\phi \otimes_R id} \\
L^{\prime} \otimes_R P \ar[rr]^{\tau(L^{\prime})} && L^{\prime} \otimes_R Q. \\
}
\end{equation}
Hence two isomorphic $(R,S)$-progenerators lead to two categorical equivalences, which are in the same isomorphism class. Thus we have Morita III.

To prove that Morita III implies Theorem \ref{Bicat}, assume that 
an invertible $R$-$S$ bimodule exists. Proposition \ref{crucial} provides us with a $(R,S)$-progenerator. Applying Morita III, we obtain a categorical equivalence $\mathfrak{M}_S \rightarrow \mathfrak{M}_R$.    

In opposite direction, starting with $\mathfrak{M}_S \simeq \mathfrak{M}_R$, we have a categorical equivalence $\mathfrak{M}_S \rightarrow \mathfrak{M}_R$. 
Applying Morita III, we obtain an $(R,S)$-progenerator. Proposition \ref{crucial} shows that this $(R,S)$-progenerator is invertible, 
so $R \biso S$. Thus we have Theorem \ref{Bicat}.
\end{subsection}
\end{section}

%% file: vNa.tex
\begin{section}{Morita theory for von Neumann algebras}
This section discusses the von Neumann algebraic analogue of the bicategorical approach to Morita theory for rings.
We will briefly repeat the necessary definitions and theory. An extensive treatment of the theory of von Neumann algebras may be found in ~\cite{B&R},~\cite{Conn1},~\cite{ev},~\cite{K&R1},~\cite{StrZs} or~\cite{sund}. 
 
\begin{subsection}{Basic definitions}
\begin{stelling}[von Neumann's Bicommutant Theorem]
Let $\mathcal{H}$ be a \\ \hbox{Hilbert space}, and let $\mathcal{B}(\mathcal{H})$ the space of all bounded linear operators on $\mathcal{H}$.
Let $\mathcal{A}$ be a $^{\star}$-subalgebra of $\mathcal{B}(\mathcal{H})$ containing $1$. Then the following conditions are equivalent:
\begin{enumerate}
\item{
$\mathcal{A}$ is closed in the $\sigma$-weak topology.}
\item{
$\mathcal{A}^{\prime\prime} = \mathcal{A}$, where the commutant $\mathcal{A}^{\prime}$ is defined by 
$$ \mathcal{A}^{\prime} = \{ x \in \mathcal{B}(\mathcal{H})\ \  |\ \  ax = xa \ \ \forall a \in \mathcal{A} \}.$$
}
\end{enumerate}
\label{bicomm}
\end{stelling}
Recall that a net $\{x_i\}_i$ converges to $x \in \mathcal{B}(\mathcal{H})$ in the $\sigma$-weak topology if 
\begin{equation}
\sum_{j = 1}^{\infty} <(x_i - x)\eta_j, \zeta_j> \rightarrow 0,
\end{equation}
for all $\eta_j,\zeta_j \in \mathcal{H}$ such that $\sum_{j = 1}^{\infty}||\eta_j||^2 < \infty$ and $\sum_{j = 1}^{\infty}||\zeta_j||^2 < \infty$. 
In fact, an even stronger result than stated in the theorem above holds. The conditions of Theorem \ref{bicomm} are also equivalent to the conditions that 
$\mathcal{A}$ is weakly closed, $\mathcal{A}$ is strongly closed, or $\mathcal{A}$ is $\sigma$-strongly closed. 
A proof of Theorem \ref{bicomm}, (as well as the definitions of the topologies mentioned above), may be found in every book on von Neumann algebras. 

\begin{definitie}
A \textdef{von Neumann algebra}\footnote{A von Neumann algebras was originally called a \emph{Ring of Operators} and later a \emph{$W^*$-algebra}; this explains our notation $\mathsf{[W^*]}$ for the bicategory of von Neumann algebras.} is a $^{\star}$-subalgebra of $\mathcal{B}(\mathcal{H})$ containing $1$ and satisfying one (and hence both) of 
the conditions of Theorem \ref{bicomm}.
\end{definitie}

Let $\m$ be a von Neumann algebra acting on a Hilbert space $\mathcal{H}$. The \emph{predual} $\m_*$, consists of all normal linear functionals on $\m$. Recall that normality is equivalent to $\sigma$-weak continuity for functionals and representations. In what follows, we will not distinguish between the two (see~\cite{B&R}, Lemma 2.4.19). 
A crucial property of von Neumann algebras is that a von Neumann algebra $\m$ is characterized by its predual in the following way:
The predual $\m_*$ is the (unique) Banach space, for which the dual is isomorphic to $\m$: $(\m_*)^* \simeq \m$. Furthermore, the $\sigma$-weak topology on $\m$ coincides with the weak-$^*$ topology on $(\m_*)^*$ as a Banach space. (See~\cite{ev} Thm. 5.11).

\begin{voorbeeld}
For a $(X,\mu)$ a measure space, consider the Hilbert space $L^2(X,\mu)$. Any $L^{\infty}$ function acts on $L^2(X,\mu)$ by multiplication. Hence $L^{\infty}(X,\mu)$ 
can be seen as a subalgebra of $\mathcal{B}(L^2(X,\mu))$; it is in fact a von Neumann algebra. The predual  $L^{\infty}(X,\mu)_*$ is given by the Banach space 
$L^1(X,\mu)$, since we have $ L^1(X,\mu)^* = L^{\infty}(X,\mu)$. 
\label{theexample}
\end{voorbeeld}
In the light of the previous example, a general von Neumann algebra may be seen as the non-commutative 
analogue of $L^{\infty}(X,\mu)$.

\begin{definitie}
Let $\m,\n$ be von Neumann algebras, $\mathcal{H}$ a Hilbert space. Suppose $\pi_l$ is a normal unital representation of $\m$ on $\mathcal{H}$ and 
$\pi_r$ is a normal unital representation of $\n^{op}$ on $\mathcal{H}$ (or equivalently, an anti-representation of $\n$) such that the actions of $\pi_l(\m)$ and $\pi_r(\n)$ commute. The triple $\bigl[\pi_l, \pi_r, \mathcal{H}\bigr]$ is called a  
\textdef{correspondence}, denoted by $\m \rightarrow \mathcal{H} \leftarrow \n$. We write $x \eta$ instead
 of $\pi_l(x)\eta$ and $\eta y$ instead of $\pi_r(y)\eta $, for $x \in \m, \eta \in \mathcal{H}, y \in \n$.
\end{definitie}

Viewing a von Neumann algebra as a ring, a correspondence may be seen as a 
bimodule. Examples of correspondences are not difficult to find. Suppose we 
have a normal, unital representation of a von Neumann algebra $\m$ on a Hilbert space $\mathcal{H}$. 
Then we immediately have a correspondence
\begin{equation}
\m \rightarrow \mathcal{H} \leftarrow (\m^{\prime})^{op}.
\end{equation}
Trivially, the actions of $\m$ and $(\m^{\prime})^{op}$ commute, and $\m^{\prime}$ acts on $\mathcal{H}$ by definition. 
Another way of constructing correspondences is through $*$-homomorphisms:
A normal $*$-homomorphism $\rho : \m \rightarrow \n$ between two von Neumann algebras gives rise to a $\n$-$\m$ correspondence. This correspondence, denoted by 
$\mathfrak{L}^2(\rho)$, is defined by 
\begin{equation} 
\mathfrak{L}^2(\rho) = \{ \xi \in \mathfrak{L}^2(\n)\ \  |\ \  \xi\rho(1) = \xi \},
\end{equation}
where $\mathfrak{L}^2(\n)$ is the standard form of $\n$, see below.
The left $\n$ representation is given by $\pi_l(n)\xi = n\xi$, and the right $\m$ representation is given by $\pi_r(m) = \xi\rho(m)$, for $ n \in \n, m \in \m, \xi \in \mathfrak{L}^2(\rho)$.  
If $\rho$ is unital, $\mathfrak{L}^2(\rho) = \mathfrak{L}^2(\n)$, otherwise $\rho(1)$ is a projection on $\mathfrak{L}^2(\n)$.
If $\n$ is properly infinite, then every $\n$-$\m$ correspondence is equivalent to an $\mathfrak{L}^2(\rho)$, for some normal 
$*$-homomorphism  $\rho : \m \rightarrow \n$ (see~\cite{Conn1}).

In what follows, we will assume that the reader is familiar with the facts and constructions below. 
All topologies mentioned above are weaker than 
the norm topology, so a von Neumann algebra is in particular a $C^*$-algebra. 
Further, if the von Neumann algebra allows a normal faithful positive linear functional we may apply the GNS-construction for von Neumann algebras (\!~\cite{sund} \S 2.2). We will also use the existence of polar decomposition.

\end{subsection}

\begin{subsection}{The standard form of a von Neumann algebra}
To state the Morita theory for von Neumann algebras, the so-called \emph{identity correspondence} or \emph{standard form} for a von Neumann algebra needs to be explained. The identity correspondence will act as the 
identity arrow in the bicategory. Note the difference with the bicategory of rings, where each ring is its own identity arrow. However, a von Neumann algebra is not an element of its own representation category (see Definition \ref{defmrep}) so we need an alternative identity arrow. 
The fact that every von Neumann algebra is isomorphic to one in standard form 
is one of the main results in Tomita-Takesaki theory.

\begin{definitie}
A von Neumann algebra $\m$ acting on a Hilbert space $\mathcal{H}$ is said to be in \textdef{standard form} if there 
exists a conjugation $J: \mathcal{H} \rightarrow \mathcal{H}$ 
such that the mapping
\begin{eqnarray}
x &\mapsto& Jx^*J,
\end{eqnarray}
defines a (complex linear) isomorphism $\m \rightarrow (\m^{\prime})^{op}$, 
which is the identity on the center \hbox{$Z(\m) = \m \cap \m^{\prime}$} of $\m$.
\end{definitie}

Some mathematicians include the existence of a selfdual cone $\mathcal{P} \subset \mathcal{H}$ with some additional properties in the definition of the standard form, see~\cite{Haag},~\cite{Tak}. 
However, the existence of such a cone follows directly from our definition (\!~\cite{Haag}, Remark 2.2 ).
 
\begin{stelling}
Every von Neumann algebra is isomorphic to one in standard form, and the standard form is unique up to unitary equivalence.
\end{stelling}

We will only give a sketch of the existence part of the proof. Uniqueness up 
to unitary equivalence is proven in~\cite{Haag}. 
Let us first assume that $\m$ has a cyclic and separating vector $\xi$ in $\mathcal{H}$, which means, respectively, 
\begin{eqnarray}
\overline{\m\xi} &=& \mathcal{H};  \\
A\xi = 0 &\Rightarrow& A = 0,\ \  A \in \m. \label{separating}
\end{eqnarray}
Recall that the vector $\xi$ being separating for $\m$ implies that $\xi$ is cyclic for $\m^{\prime}$. So equation \eqref{separating} is equivalent to 
\begin{equation}
\overline{\m^{\prime}\xi} = \mathcal{H}.
\end{equation}
Note that not every von Neumann algebra has such a cyclic and separating vector. However, if the von Neumann algebra admits a faithful normal state, we are able to construct one, using the GNS representation\footnote{For example, all von Neumann algebras that are $\sigma$-finite allow a cyclic and separating vector, see~\cite{B&R}, Prop. 2.5.6. This includes all von Neumann algebras that allow a faithful representation on a separable Hilbert space.}.  
Suppose we have a state $\phi$ on a 
von Neumann algebra $\m$. Recall that a state is a positive, normalized linear 
functional on $\m$. A functional is normalized if $\phi(1) = 1$. 
It is normal when 
\begin{equation}
\phi(\sup_i x_i) = \sup_i \phi(x_i),
\label{normal}
\end{equation}
for any increasing net $\{x_i\}_i$ in the positive cone $\m_+$. Applying the GNS-construction to $(\m,\phi)$, we obtain a Hilbert space $\mathcal{H}_{\phi}$, a normal representation $\pi_{\phi}$ on $\mathcal{H}_{\phi}$, and a cyclic and separating vector $\xi_{\phi}$.    

For now, let $\m$ be a von Neumann algebra acting on an Hilbert space 
$\mathcal{H}$, and let $\xi$ be a cyclic and separating vector for $\m$. 
Define unbounded anti-linear operators 
\begin{equation}
Sx\xi = x^* \xi,
\end{equation} 
for $x \in \m$ on the dense domain $\m\xi$, and
\begin{equation}
Fy'\xi = {y'}^* \xi,
\end{equation}
for $y' \in \m^{\prime}$ on the dense domain $\m^{\prime}\xi$. A trivial computation shows that
\begin{equation}
<Sx \xi,y' \xi> = <Fy' \xi,x \xi>,
\end{equation}
for $x \in \m, y' \in \m^{\prime}$,
so $F = S^*$ as antilinear operators. Hence both $S$ and $F$ have a densely defined linear conjugate. This means $S$ and 
$F$ are closable and we can apply the polar decomposition to the closures, again 
denoted by $S$ and $F$.
Let
\begin{equation}
S = J \Delta^{1/2},
\end{equation}
be the polar decomposition. We have the identities
\begin{eqnarray}
S &=& \Delta^{-1/2}J \nonumber \\
J &=& J^*\  = \ J^{-1} \nonumber \\
S m^* \eta &=& m S\eta, \nonumber \\
\Delta^{1/2}m &=& m\Delta^{1/2} ,
\label{rekenregels}
\end{eqnarray}
$\eta \in \mathcal{H}, m \in Z(\m)$.
The first two equalities follow easily from the polar decomposition, proof of the third can be found in~\cite{StrZs} p. 273. A proof of the fourth 
property may be found in Lemma \ref{lemmaabove} below.
Since $S, \Delta$ are densely defined, the partial 
isometry $J$ is an anti-linear isometry on $\mathcal{H}$.
 Now $J$ is a self-adjoint, anti-unitary 
operator, and $\Delta$ is a positive self-adjoint (unbounded) operator. 
A rather lengthy argument shows that 
\begin{equation}
\m  \simeq (\m^{\prime})^{op}, 
\end{equation}
through
\begin{eqnarray}
x \mapsto Jx^*J, \label{zm}
\end{eqnarray}
so that 
\begin{equation}
J\m J = (\m^{\prime})^{op},
\end{equation}
as algebras.
See~\cite{sund} for
the proof where $S$ is assumed bounded and~\cite{K&R} for 
the general case. Using the identities \eqref{rekenregels}, 
it is easy to show that the isomorphism \eqref{zm} indeed acts like 
the identity on the center. Let $x \in Z(\m), \eta \in \mathcal{H}$, 
\begin{eqnarray}
Jx^*J\eta &=& Jx^*\Delta^{1/2}S\eta = J\Delta^{1/2}x^*S\eta \nonumber \\
&=& Sx^*S\eta = S^2 x \eta = x \eta.
\end{eqnarray}

In the case that we do not have a faithful normal 
state, we can use a similar construction 
using a weight:
\begin{definitie}
A \textdef{weight} on a von Neumann algebra $\m$ is an additive map $\phi$ from the positive cone $\m_+$ into the extended reals $\overline{\mathbb{R}_+}= [0,\infty]$. 
The map $\phi$ is called \textdef{semifinite} if
$$ \mathfrak{D}_{\phi} = \{ x \in \m_+ \ \ | \ \ \phi(x) < \infty\} $$
generates $\m$ as a von Neumann algebra, 
and \textdef{faithful} if $\phi(x) = 0, x \in \m_+$ implies $x = 0$.
Let $\{x_i\}_i$ be a bounded increasing net in $\m_+$. Then $\phi$ is called
\textdef{normal} if \eqref{normal} holds. 
\end{definitie}

We would like to apply a GNS-like construction to obtain a Hilbert space and 
a representation. However, this construction will not provide us with a 
cyclic and separating vector, since the identity operator is not an element of $\mathfrak{D}_{\phi}$ if 
the weight is not finite.

Let $\phi$ be an arbitrary faithful semifinite normal weight on $\m$. Define 
\begin{equation}
\m_{\phi} = \{ x \in \m \ \ | \ \ \phi(x^*x) < \infty\}.
\end{equation}

Note that $\m_{\phi}$ is a left ideal 
of $\m$, since  $\phi$ is an additive positive mapping and one has 
\begin{eqnarray} 
(x + y)^*(x + y) &\leq& 2(x^*x + y^*y), \\
(mx)^*(mx) &\leq& ||m||_{\m}^2\ x^*x \label{deze!},
\end{eqnarray}
for $x,y \in \m_{\phi}, m \in \m$.
Hence we have an action of $\m$ by left multiplication on the space 
$\m_{\phi}$.  Next, we define an inner product on $\m_{\phi}$ by 
\begin{equation}
(x,y) = \phi(x^*y),  
\label{inprod}
\end{equation}
for $x,y \in \m_{\phi}$. 
It is easy to check that the inner product \eqref{inprod} is well defined and 
positive definite. The completion of $\m_{\phi}$ will be denoted 
by $\mathcal{H}_{\phi}$. Using \eqref{deze!} 
we have for $y \in \m, x \in \m_{\phi}$,
\begin{equation}
{||yx||}_{\mathcal{H}_{\phi}}^2 \leq ||y||_{\m}^2\  \phi(x^*x),
\end{equation}
which implies that the left multiplication of $\m$ on $\m_{\phi}$ is 
bounded, so we can extend it to a bounded operator on $\mathcal{H}_{\phi}$. 
Hence 
\begin{eqnarray}
\pi_{\phi} : \m &\rightarrow& \mathcal{B}(\mathcal{H}_{\phi}); \nonumber \\
\pi_{\phi}(x)\eta &=& x\eta, 
\end{eqnarray}
is a representation of $\m$ on $\mathcal{H}_{\phi}$, and, by construction, this representation is faithful.

\begin{opmerking}
Note that it is not necessary for $\phi$ to be faithful. If we define 
\begin{equation}
\mathcal{N}_{\phi} = \{ x \in \m \ \ | \ \ \phi(x^*x) = 0 \},
\end{equation}
which is a left ideal of $\m$ as well, we can define an inner product on the 
quotient space $\m_{\phi}/\mathcal{N}_{\phi}$. Similarly to the construction above, we obtain a Hilbert space which is the completion of the 
quotient in this inner product and a 
faithful representation.
\end{opmerking}

It can be shown that $\m_{\phi} \cap {\m_{\phi}}^*$ with the structure of a 
$^*$-algebra induced by $\m$, and the scalar product induced by $\mathcal{H}_{\phi}$, 
is a left Hilbert algebra $\mathfrak{A} \subset \mathcal{H}_{\phi}$;
the associated von Neumann 
algebra $\mathcal{L}(\mathfrak{A})$ is isomorphic to $\pi_{\phi}(\m)$. 
A left Hilbert algebra admits a preclosed anti-linear operator 
\begin{eqnarray}
S : \mathfrak{A} &\rightarrow& \mathfrak{A} \nonumber \\
\eta &\mapsto& \eta^*,
\label{moeizaam}
\end{eqnarray}
given by the involution on $\mathfrak{A}$. As above, the polar decomposition 
$S = J\Delta^{1/2}$ provides an anti-unitary conjugation $J$ and a positive self-adjoint (unbounded) operator $\Delta$. The operators $J$
 and $\Delta$ are called the \emph{modular conjugation} and 
\emph{modular operator} associated with $\m$, respectively.
The anti-unitary conjugation $J$ in $\mathcal{H}_{\phi}$,
 defines an anti-isomorphism  
\begin{equation}
x \mapsto Jx^*J, \hspace{20pt} x \in  \mathcal{L}(\mathfrak{A}),
\end{equation}
from
$\mathcal{L}(\mathfrak{A})$ to  ${\mathcal{L}(\mathfrak{A})}^{\prime}$, 
which acts like the identity on the center. This proves that $\m$ is in standard 
form, since we can identify $\m$ with $\pi_{\phi}(\m)$. 
An exact definition of a left Hilbert algebra and a proof of the statements above can be found in~\cite{StrZs}. 
A brief discussion of the standard form, starting with left Hilbert 
algebras can be found in~\cite{Tak}. 

It is left to show that every von Neumann algebra admits a faithful, normal, 
semifinite weight (see~\cite{K&R},~\cite{StrZs}). Let ${\{\phi_i\}}_i$ be a maximal family of normal forms (i.e. normal linear functionals) on 
$\m$, whose supports are mutually orthogonal. 
The formula
\begin{equation}
\phi(a) = \sum_i \phi_i(a), \hspace{20pt} a \in \m_+, 
\end{equation}
yields a faithful, semifinite, normal weight on $\m_+$.
By definition, $\phi$ is a weight. Faithfulness follows directly from the facts 
that all $\phi_i$ are positive and that it is a maximal family. 
Let ${\{x_j\}}_j$ be a bounded increasing 
net in $\m_+$.
Then
\begin{eqnarray}
\sup_j \phi(x_j) &=&  \sup_j \sum_i \phi_i(x_j) = \sum_i \sup_j \phi_i(x_j)
\nonumber \\ &=& \sum_i \phi(\sup_j x_j) = \phi(\sup_j x_j), 
\end{eqnarray}
so $\phi$ is normal.   
Semifiniteness of $\phi$ follows from the fact that 
\begin{equation}
1 = \bigvee   \{ p \in \mathcal{P}(\m)\ \  | \ \ \phi(p) < \infty \},
\label{sundeq}
\end{equation}
where $\mathcal{P}(\m)$ denotes the set of all projections in $\m$, and $\bigvee_i p_i$ is the projection on the subspace $\overline{\sum_i p_i\mathcal{H}}$. 
which holds since the supports of ${\{\phi_i\}}_i$ are mutually orthogonal and the family ${\{\phi_i\}}_i$
 is maximal. Equation \eqref{sundeq} is equivalent to the fact that $\phi$ is semifinite (\!~\cite{sund} p. 57).
\end{subsection}

\begin{subsection}{The identity correspondence}

\begin{definitie}
Let $\m$ be a von Neumann algebra acting on a Hilbert space $\mathcal{H}$ 
which is in standard form. Then the \textdef{identity correspondence} of 
$\m$ is the Hilbert space $\mathcal{H}$, the left representation $\pi_l$ given 
by multiplication and the right representation $\pi_r$ defined via the anti-unitary 
conjugation $J$:
\begin{eqnarray}
\pi_l : \m &\longrightarrow& \mathcal{B}(\mathcal{H}); \nonumber \\
\pi_l(x) &\longmapsto& (\eta \mapsto x\eta),
\end{eqnarray}
and
\begin{eqnarray}
\pi_r : \m  &\longrightarrow& \mathcal{B}(\mathcal{H}); \nonumber \\
\pi_r(x) &\longmapsto& (\eta \mapsto Jx^*J\eta). \label{rechtsrep}
\end{eqnarray}
\end{definitie}


The identity correspondence is denoted by 
$$\m \rightarrow \mathfrak{L}^2(\m) \leftarrow \m.$$
The notation $\mathfrak{L}^2(\m)$ is chosen analogously to measure theory; recall that it is considered to be the commutative version of von Neumann algebras, see Example \ref{theexample}.
\end{subsection}

\begin{subsection}{The relative tensor product or Connes fusion}
This section handles the relative tensor product of two matching correspondences. In the construction of the bicategory $\mathsf{[W^*]}$, classes of correspondences will form the horizontal arrows in the bicategory, analogously to classes of bimodules in the bicategory $\mathsf{[Rings]}$. The relative tensor product will play the role 
of the composition functor, analogously to the tensor product of bimodules in the case of rings.  

Let $\m \rightarrow \mathcal{H} \leftarrow \n$ and $\n \rightarrow \mathcal{K} \leftarrow \p$ be two correspondences. 
To obtain the so-called \emph{relative 
tensor product} of $\mathcal{H}$ and $\mathcal{K}$ we follow~\cite{sau}. However, this construction is not symmetric in $\mathcal{H}$ and $\mathcal{K}$. 
Wassermann~\cite{wass} has stated a symmetric construction in the case that the von Neumann algebra allows a (cyclic and separating) vacuum vector. The construction below 
is applicable to general von Neumann algebras. Consider a normal, faithful, semifinite weight $\psi$ on $\n$. Recall the 'GNS'-like construction from the previous section, 
which provides us with a Hilbert space $\mathcal{H}_{\psi}$, a dense subset $$\n_{\psi} = \{ y \in \n \ \ | \ \ \psi(y^*y) < \infty\},$$ and a representation $\pi_{\psi}$.
Let $\Lambda_{\psi} : \n_{\psi} \rightarrow \mathfrak{L}^2(\n)$ denote the canonical inclusion map. 
We will often abuse this notation, and write $y \in \n_{\psi}\subset \mathfrak{L}^2(\m)$ 
when we consider $y$ as $\Lambda(y)$, an element of the Hilbert space $\mathfrak{L}^2(\n)$.

Define the subset $\widetilde{\mathcal{H}} \subset \mathcal{H}$ of \textdef{$\psi$-bounded} vectors as the set of $\eta \in \mathcal{H}$ for which 
\begin{equation} 
{||\eta y||}_{\mathcal{H}} \leq C_{\eta} {||y||}_{\mathcal{H}_{\psi}}, \hspace{20pt} \forall y \in \n_{\psi},
\end{equation}
holds.
Note that ${||y||}_{\mathcal{H}_{\psi}}$ is finite by definition of $\n_{\psi}$ and of the norm on $\mathcal{H}_{\psi}$. 
Equivalently, one can define the subspace $\widetilde{\mathcal{H}}$ as the set of $\eta \in \mathcal{H}$ for which the operator 
\begin{equation}
R^{\psi}_\eta : \mathcal{H}_{\psi} \rightarrow \mathcal{H},
\end{equation}
defined on the dense subspace $\n_{\psi}$ by
\begin{equation}
R^{\psi}_{\eta}(J\Lambda(y^*)) = \pi_r(y)\eta  \hspace{20pt} \forall y \in \n_{\psi},
\label{defJ}
\end{equation}
is bounded. The operator $J$ in \eqref{defJ} is the modular conjugation associated with $\n$. Observe that $J$ originates from an operator 
$S : \n_{\psi} \cap \n_{\psi}^* \rightarrow \mathcal{H}_{\psi}$ (see \eqref{moeizaam}), since we can not assume that 
$\n$ allows a cyclic and separating vector.

Recall that a construction like above yields the standard form of $\n$, so we may also write $\mathfrak{L}^2(\n)$ instead of $\mathcal{H}_{\psi}$ and 
\begin{equation}
R_{\eta} : \mathfrak{L}^2(\n) \rightarrow \mathcal{H},
\end{equation}
where $R_{\eta}$ is defined on a dense subspace of $\mathfrak{L}^2(\n)$ by the right representation of $\n$ on $\mathcal{H}$. 
In the light of the previous remark, it seems reasonable that 
the construction of the relative tensor product does not depend on the weight $\psi$ (up to unitary equivalence). 
This is indeed the case, a proof can be found in~\cite{sau}. 
\begin{lemma} 
Let $\m \rightarrow \mathcal{H} \leftarrow \n$ be a correspondence, and define   
$\widetilde{\mathcal{H}}$  with the use of an arbitrary normal faithful semifinite weight $\psi$. Let $\Delta$ be the modular operator associated with $\n$. 
The following properties hold.
\begin{enumerate}
\item{$\widetilde{\mathcal{H}}$ is stable under the actions of $\m$ and $\n$.} \label{1e}
\item{$\widetilde{\mathcal{H}}$ is a dense subspace of $\mathcal{H}$.} \label{2e}
\item{For $\eta_1,\eta_2 \in \widetilde{\mathcal{H}}$, $R_{\eta_1}^{*\phantom{1}} R_{\eta_2}^{\phantom{1}} \in \n$ holds, where $\n$ is identified with the left representation on $\mathfrak{L}^2(\n)$.} \label{3e}
\item{For $\eta_1, \eta_2 \in \widetilde{\mathcal{H}}, B \in (\n^{op})^{\prime} \subset \mathcal{B}(\mathcal{H})$, 
we have that $R_{B\eta_1}^{*\phantom{B1}}R_{B\eta_2}^{\phantom{B\eta}}$ equals $R_{\eta_1}^{*\phantom{1}}B^*BR_{\eta_2}^{\phantom{1}}$ as 
operators on $\mathfrak{L}^2(\n)$.}\label{4e}
\item{For $\eta_1, \eta_2 \in \widetilde{\mathcal{H}}, A \in \n$, we have that $\Delta^{-1/2}A^*\Delta^{1/2}R_{\eta_1}^{*\phantom{1}} R_{\eta_2}^{\phantom{1}}$ 
equals $R_{\eta_1 A}^{*\phantom{A}} R_{\eta_2}^{\phantom{!}}$ as operators on $\mathfrak{L}^2(\n)$. Analogously, $R_{\eta_1}^{*\phantom{1}}R_{\eta_2}^{\phantom{e}}\Delta^{1/2}
A\Delta^{-1/2}$ equals $R_{\eta_1}^{^*\phantom{1}}R_{\eta_2 A}^{\phantom{asd}}$ as operators on $\mathfrak{L}^2(\n)$.} \label{help}

\end{enumerate}
\begin{bewijs}
Proof of \eqref{1e}:
First we will show that $\widetilde{\mathcal{H}}$ is stable under $\n$.
For \\$A \in \n, \eta \in \widetilde{\mathcal{H}}, 
y \in \n_{\psi}$ we have
\begin{eqnarray}
{||R^{\psi}_{\eta A}(Jy^*)||}_{\mathcal{H}} &=&  {||\pi_r(y)(\eta A)||}_{\mathcal{H}} =  {||(\eta A)y||}_{\mathcal{H}} 
\leq  C_{\eta} {||Ay||}_{\mathcal{H}_{\psi}} \nonumber \\ &\leq& 
C_{\eta}  {||A||}_{\mathcal{H}} {||y||}_{\mathcal{H}_{\psi}} = C_{\eta,A}
 {||y||}_{\mathcal{H}_{\psi}},
\end{eqnarray}
where we used that $\n_{\psi}$ is an ideal in $\n$. The  
second inequality holds since
\begin{equation}
||Ay||^2_{\mathcal{H}_{\psi}} = \psi\bigl((Ay)^*(Ay)\bigr) \leq ||A||^2_{\mathcal{H}} \psi(y^*y).
\end{equation}
Since $\m \subset (\n^{\prime})^{op}$, it is sufficient to show that $\widetilde{\mathcal{H}}$ is stable under $\n^{\prime}$.
Let \\ $B \in \n^{\prime}, \eta \in \widetilde{\mathcal{H}}, 
y \in \n_{\psi}$.
We have
\begin{eqnarray}
{||R^{\psi}_{B\eta}(Jy^*)||}_{\mathcal{H}} &=&  {||\pi_r(y)(B\eta)||}_{\mathcal{H}} = {||(B\eta) y||}_{\mathcal{H}} = 
{||(\eta y) B||}_{\mathcal{H}} \nonumber \\
 &\leq& {||\eta y||}_{\mathcal{H}} {||B||}_{\mathcal{H}} \leq C_{\eta B} {||y||}_{\mathcal{H}_{\psi}}. 
\end{eqnarray}

Proof of \eqref{2e}: 
This statement is due to Connes, see~\cite{Conn2}. For $n \in \n$, let $\pi(n)$ denote the corresponding operator in $\mathcal{H}$ and let $k \in \n$ denote 
the central projection corresponding to the kernel of the representation $\pi$. Now $\pi$ is an isomorphism of $\n_{(1-k)}$ to $\pi(\n)$. Let $\psi_{(1-k)}$ denote the 
restriction of $\psi$ to $\n_{(1-k)}$. Applying~\cite{ped}, Cor. 4.6.12, 
there exists a family $\{\xi_{\alpha}\}_{\alpha \in I}$ of vectors in $\mathcal{H}$  such that 
\begin{equation}
\psi_{(1-k)}(n) = \sum_{\alpha}<\pi(n)\xi_{\alpha},\xi_{\alpha}>
\label{phi=inprod}
\end{equation}     
for all $n \in (\n_{(1-k)})_+$, since $\psi_{(1-k)}$ is normal.
The family $\{\xi_{\alpha}\}_{\alpha}$ consists of $\psi$-bounded vectors, since 
\begin{equation}
||\pi(n)\xi_{\alpha}||^2 \leq C_{\xi_{\alpha}}\psi(n^*n),
\end{equation}
for all $n \in \n$. 
Let $E$ be the orthogonal projection of $\mathcal{H}$ on the closure of $\widetilde{\mathcal{H}}$. We will show that $E=1$ on $\n_{(1-k)}$, 
which finishes the proof. The projection $E$ commutes with $\pi(\n)^{\prime}$, so $E \in \pi(\n)^{\prime\prime} = \pi(\n)$. Therefore, $E$
 must be of the form $E = \pi(e)$ for some $e \in \n_{(1-k)}$. 
It now follows from \eqref{phi=inprod} that, if $e \neq 1-k$ then $\psi((1-k) - e) > 0$, and 
$\exists \ \xi_{\alpha}$ such that $\pi(1- k - e)\xi_{\alpha} > 0$
which contradicts $E\xi_{\alpha} = \xi_{\alpha}$. 

Proof of \eqref{3e}:
For $\eta_1,\eta_2 \in \widetilde{\mathcal{H}}$, we have 
\begin{eqnarray}
R_{\eta_1}^{*\phantom{1}}R_{\eta_2}^{\phantom{*}} \in Hom_{\n^{op}}(\mathfrak{L}^2(\n), \mathfrak{L}^2(\n)),
\end{eqnarray}
since $R_{\eta_1}^{*\phantom{1}}R_{\eta_2}^{\phantom{1}} : \mathfrak{L}^2(\n) \rightarrow \mathfrak{L}^2(\n)$ by definition, and 
$R_{\eta}$ intertwines the right $\n$-action. The computation 
\begin{eqnarray}
Hom_{\n^{op}}(\mathfrak{L}^2(\n), \mathfrak{L}^2(\n)) =  (\n^{op})^{\prime} = \n^{\prime\prime} = \n
\end{eqnarray}
finishes the proof.

Proof of \eqref{4e}:
 Let $\eta_1,\eta_2 \in \widetilde{\mathcal{H}}, B \in (\n^{op})^{\prime}$. 
Let $x,y \in \n_{\psi}$, and let $( , )$ denote the inner product on $\mathfrak{L}^2(\n)$.
Then
\begin{eqnarray}
(R_{B\eta_1}^{*\phantom{1}}R_{B\eta_2}^{\phantom{B\eta}}Jx,Jy) &=& (R_{B\eta_1}^{*\phantom{B\eta}}B\eta_2 x^*,Jy) \nonumber \\
&=& (R_{B\eta_1}^{*\phantom{B\eta}}BR_{\eta_2}^{\phantom{1}}Jx,Jy) \nonumber \\
&=& (Jx, R_{\eta_2}^{*\phantom{1}}B^*R_{B\eta_1}Jy) \nonumber \\
&=& (Jx, R_{\eta_2}^{*\phantom{1}}B^*B\eta_1 y^*) \nonumber \\
&=& (Jx, R_{\eta_2}^{*\phantom{1}}B^*B R_{\eta_1} Jy) \nonumber \\
&=& (R_{\eta_1}^{*\phantom{2}}B^*B R_{\eta_2}^{\phantom{1}}Jx,Jy).
\end{eqnarray}
Remark that we have implicitly used that $\widetilde{\mathcal{H}}$ is closed under the $(\n^{op})^{\prime}$-action; c.f. the proof of \eqref{2e} above.
 Now $R_{B\eta_1}^{*\phantom{B1}}R_{B\eta_2}^{\phantom{B\eta}}$ equals $R_{\eta_1}^{*\phantom{1}}B^*BR_{\eta_2}^{\phantom{1}}$ as 
operators on $\mathfrak{L}^2(\n)$, since $\n_{\psi}$ is dense in 
$\mathfrak{L}^2(\n)$.

Proof of \eqref{help}: Let $\eta_1,\eta_2 \in \widetilde{\mathcal{H}}, A \in \n$.
 Let $x,y \in \n_{\psi}\cap \n_{\psi}^*$. Then
\begin{eqnarray}
(R^{*\phantom{A}}_{\eta_1 A}R_{\eta_2}^{\phantom{**}}Jx,Jy) &=& (Jx,R_{\eta_2}^{*\phantom{@}}R_{\eta_1 A}^{\phantom{**}}Jy) \nonumber \\
&=& (Jx,R_{\eta_2}^{*\phantom{@}}\eta_1 A y^*) \nonumber \\
&=& (Jx,R_{\eta_2}^{*\phantom{@}}R_{\eta_1}^{\phantom{**}}J(Ay^*)^*) \nonumber \\
&=& (Jx,R_{\eta_2}^{*\phantom{@}}R_{\eta_1}^{\phantom{**}}\Delta^{1/2}S(Ay^*)^*) \nonumber \\
&=& (Jx,R_{\eta_2}^{*\phantom{@}}R_{\eta_1}^{\phantom{**}}\Delta^{1/2}Ay^*) \nonumber \\
&=& (Jx,R_{\eta_2}^{*\phantom{@}}R_{\eta_1}^{\phantom{**}}\Delta^{1/2}ASy) \nonumber \\
&=& (Jx,R_{\eta_2}^{*\phantom{@}}R_{\eta_1}^{\phantom{**}}\Delta^{1/2}A\Delta^{-1/2}Jy) \nonumber \\
&=& (\Delta^{-1/2}A^*\Delta^{1/2}R_{\eta_1}^{*\phantom{a}}R_{\eta_2}^{\phantom{**}}Jx,Jy).
\label{deltascore}
\end{eqnarray}
Note that we take $x,y \in \n_{\psi}\cap \n_{\psi}^*$, otherwise the expression $Sy$ would not be defined. Now $R^{*\phantom{A}}_{\eta_1 A}R_{\eta_2}$ equals 
$\Delta^{-1/2}A^*\Delta^{1/2}R_{\eta_1}^{*\phantom{a}}R_{\eta_2}^{\phantom{**}}$ as operators on $\mathfrak{L}^2(\n)$ since $\n_{\psi}\cap \n_{\psi}^*$ is dense in $\mathcal{H}_{\psi} = 
\mathfrak{L}^2(\n)$. Analogously, we have 
\begin{eqnarray}
(R_{\eta_1}^{*\phantom{1}}R_{\eta_2 A}^{\phantom{Kg}}Jx,Jy) &=& (Jx,R_{\eta_2 A}^{*\phantom{@D}}R_{\eta_1}^{\phantom{**}}Jy) \nonumber \\
&=& (Jx, \Delta^{-1/2}A^*\Delta^{1/2}R_{\eta_2}^{*\phantom{a}}R_{\eta_1}^{\phantom{**}}Jy) \nonumber \\
&=& (R_{\eta_1}^{*\phantom{*}}R_{\eta_2}^{\phantom{*a}}\Delta^{1/2}A\Delta^{-1/2}Jx,Jy),
\label{deltascoresagain}
\end{eqnarray}
where we used \eqref{deltascore}. Again, a density argument shows that $R_{\eta_1}^{*\phantom{1}}R_{\eta_2 A}^{\phantom{Kg}}$ 
 equals $R_{\eta_1}^{*\phantom{*}}R_{\eta_2}^{\phantom{*a}}\Delta^{1/2}A\Delta^{-1/2}$ as operators on $\mathfrak{L}^2(\n)$.
\label{5e}
\end{bewijs}
\label{c4}
\end{lemma}

We now resume the construction of the relative tensor product.
Given the correspondences  $\m \rightarrow \mathcal{H} \leftarrow \n$ and $\n \rightarrow \mathcal{K} \leftarrow \p$, first construct the 
algebraic tensor product $\widetilde{\mathcal{H}} \otimes_{\mathbb{C}} \mathcal{K}$. On this space, define a sesquilinear form, by extension of 
\begin{equation}
\bigl(\eta_1 \otimes \zeta_1,\eta_2 \otimes \zeta_2\bigr)_0 \ := \ 
<\zeta_1,R_{\eta_1}^{*\phantom{2}}R_{\eta_2}^{\phantom{1}}\zeta_2>_{\mathcal{K}},
\end{equation}
where $<\ ,\ >_{\mathcal{K}}$ is the inner product on the Hilbert space $\mathcal{K}$. This makes $( \ \ ,\ \ )_0$ a pre-inner 
product on $\widetilde{\mathcal{H}} \otimes_{\mathbb{C}} \mathcal{K}$. Note that, since $R_{\eta_1}^{*\phantom{1}}R_{\eta_2}^{\phantom{1}} \in \n$ by Lemma \ref{c4}.\ref{3e}, the second argument of the 
inner product is indeed an element of $\mathcal{K}$. 
Let $\mathcal{N}$ be the null space of the this pre-inner product. 
\begin{equation}
\mathcal{N} = \{ \eta \otimes \zeta \in \widetilde{\mathcal{H}} \otimes_{\mathbb{C}} \mathcal{K}
\ | \ \bigl(\eta \otimes \zeta, \eta \otimes \zeta\bigr)_0 = 0\}.
\end{equation}
Then $(\ ,\ )_0$ becomes an inner product on $\widetilde{\mathcal{H}} \otimes_{\mathbb{C}} \mathcal{K} / \mathcal{N}$, denoted by $<\ \ ,\ \ >_0$. 
Completion in this inner product forms a Hilbert 
space denoted by $\mathcal{H} \sqtimes_{\n} \mathcal{K}$, the \emph{relative tensor product} or \emph{Connes fusion} of  $\mathcal{H}$ and $\mathcal{K}$.

This Hilbert space $\mathcal{H} \sqtimes_{\n} \mathcal{K}$ is an $\m \rightarrow \mathcal{H} \sqtimes_{\n} \mathcal{K} \leftarrow \p$ correspondence, so that we may regard 
the above construction as the fusion of correspondences rather than merely Hilbert spaces. Namely, 
using the fact that $\mathcal{H}$ carries a left representation $\pi_l(\m)$ and $\mathcal{K}$ carries a right representation $\pi_r(\p)$ we can define 
\begin{equation}
\widetilde{\pi_l} = \pi_l \otimes id_{\mathcal{K}},
\end{equation}
and
\begin{equation}
\widetilde{\pi_r} = id_{\mathcal{H}} \otimes \pi_r,
\end{equation}
on $\mathcal{H} \otimes_{\mathbb{C}} \mathcal{K}$.
Since $\widetilde{\mathcal{H}}$ is stable under the action of $\m$, the image of $\m$ under the left representation is still in $\widetilde{\mathcal{H}} 
\otimes_{\mathbb{C}} 
\mathcal{K}$. (Trivially, $\mathcal{K}$ is stable under the action of $\p$.) It remains to be shown that the null space $\mathcal{N}$ is stable under the actions of 
$\m$ and $\p$. Let $B \in \m \subseteq (\n^{op})^{\prime}, C \in \p \subseteq (\n^{op})^{\prime}$ 
and \nopagebreak $\sum_i(\eta_i \otimes \zeta_i), \sum_j(\eta_j \otimes \zeta_j) \in \mathcal{N}$. Using Lemma \ref{c4}.\ref{4e}, we obtain  
\begin{eqnarray}
\bigl(B\sum_i(\eta_i \otimes \zeta_i), B\sum_j(\eta_j \otimes \zeta_j)\bigr)_0 &=& \bigl(\sum_i(B(\eta_i) \otimes \zeta_i), \sum_j(B(\eta_j) \otimes \zeta_j)\bigr)_0 \nonumber \\
&=& \sum_{i,j}<\zeta_i, R_{B\eta_i}^{*\phantom{B1}}R_{B\eta_j}^{\phantom{BI}}\zeta_j>_{\mathcal{K}} \nonumber \\
&=& \sum_{i,j}<\zeta_i, R_{\eta_i}^{*\phantom{1}}B^*BR_{\eta_j}^{\phantom{!}}\zeta_j>_{\mathcal{K}} \nonumber \\
&\leq& ||B||^2  \sum_{i,j}<\zeta_i, R_{\eta_i}^{*\phantom{1}}R_{\eta_j}^{\phantom{!}}\zeta_j>_{\mathcal{K}} \nonumber \\
&=& ||B||^2 \sum_{i,j}\bigl(\eta_i \otimes \zeta_i, \eta_j \otimes \zeta_j\bigr)_0 \nonumber \\
&=& 0.
\end{eqnarray}
We have used the inequality $A^*B^*BA \leq ||B||^2A^*A$ .
Hence the null space is stable under the action of $\m$. Further, using the Cauchy-Schwarz inequality on the pre-inner product $(\ \ ,\ \ )_0$ we have 
\begin{eqnarray}
&&|\Bigl(\sum_i(\eta_i \otimes \zeta_i)C,\sum_j(\eta_j \otimes \zeta_j)C\Bigr)_0| = |\Bigl(\sum_i(\eta_i \otimes (\zeta_iC)),\sum_j(\eta_j \otimes (\zeta_jC))\Bigr)_0| \nonumber \\
&& \ \ \ \  = |\sum_{i,j}<\zeta_iC, R_{\eta_i}^{*\phantom{1}}R_{\eta_j}^{\phantom{i}}\zeta_jC>_{\mathcal{K}}| \nonumber \\
&& \ \ \ \ = |\sum_{i,j}<\zeta_iC^*C, R_{\eta_i}^{*\phantom{!}}R_{\eta_j}^{\phantom{i}}\zeta_j>_{\mathcal{K}} |\nonumber \\
&& \ \ \ \ = | \Bigl(\sum_i(\eta_i \otimes (\zeta_i C^*C)), \sum_j(\eta_j \otimes \zeta_j)\Bigr)_0 |\nonumber \\
&& \ \ \ \ \leq \sum_i{\Bigl((\eta_i \otimes \zeta_iC^*C),(\eta_i \otimes \zeta_iC^*C)\Bigr)_0}^{1/2} \cdot  
\sum_j{\Bigl((\eta_j \otimes \zeta_j),(\eta_j \otimes \zeta_j)\Bigr)_0}^{1/2} \nonumber \\
&& \ \ \ \ = 0,
\end{eqnarray}
so the null space $\mathcal{N}$ is stable under $\p$ as well.

Combining several statements of~\cite{B&R}, Section 2.5.3, and~\cite{P&T} will prove the following 
lemma. We will need these rather technical results in the construction of the 
bicategory in the following section. The lemma asserts that the operations we will apply are well-defined. See also~\cite{sau}, where the second result is used without an explicit proof.  
\begin{lemma}
Let $\m$ be a von Neumann algebra, and $\phi$ an arbitrary normal faithful semifinite weight. Let $x \in \m$. Let $\tau_t(x) = \Delta^{it}x\Delta^{-it}$, $x \in \m$, be the modular automorphism group of $\m$. Define
\begin{equation}
x_n = \sqrt{n / \pi} \int dt \ e^{-nt^2} \tau_t(x) \hspace{10pt} n = 1,2,...
\end{equation}
and 
\begin{equation} 
\m_0 = \{ x_n \ \ | \ \ x \in \m_{\phi}\}.
\end{equation}
Then the following hold: 
\begin{enumerate}
\item{$\m_0 \subset \m_{\phi}$ as a $\sigma$-weakly dense subspace. Moreover, $\m_0 \subset \m_{\phi}$ as a norm dense subspace.}
\item{For $y \in \m_0$, we have $\Delta^{1/2}y\Delta^{-1/2} \in \m_{\phi}$. } 
\item{For $x \in Z(\m)$, we have $\Delta^{1/2}x\Delta^{-1/2} = x$.} 
\end{enumerate}
\begin{bewijs}
First, we would like to show that if $x \in \m_{\phi}$, then $x_n \in \m_{\phi}$. Consider $\phi(x_n^*x_n^{\phantom{p}})$. Using the fact that the integral is a Bochner integral 
we obtain 
\begin{eqnarray}
|\phi(x_n^*x_n^{\phantom{L}})| &=&  |\frac{\pi}{n} \phi(\int dt\ e^{-nt^2}\tau_t(x^*) \int ds\ e^{-ns^2}\tau_s(x))| \nonumber \\
&\leq& \frac{\pi}{n} \int dt\ ds\ e^{-nt^2} e^{-ns^2} |\phi(\tau_t(x^*)\tau_s(x))|, 
\label{abc}
\end{eqnarray}
see,~\cite{Yos}, Ch.V.
Further,
\begin{equation}
|\phi(\tau_t(x^*)\tau_s(x))| \leq \phi(\tau_t(x^*x))^{1/2} \cdot \phi(\tau_s(x^*x))^{1/2}, 
\label{def}
\end{equation}
since $\phi$ is positive semidefinite, and \begin{equation}
\tau_t(x^*)\tau_t(x) = \Delta^{it}x^*\Delta^{-it}\Delta^{it}x\Delta^{-it} = \tau_t(x^*x).
\end{equation}
It is a well-known property of the weight $\phi$ that it is invariant 
under the action of the modular automorphism group (\!~\cite{Tak}, p.17), so we have 
\begin{equation}
\phi(\tau_t(x^*x)) = \phi(x^*x).
\label{ghi}
\end{equation}
Combining \eqref{abc}, \eqref{def} and \eqref{ghi} shows that $x_n \in \m_{\phi}$.
From~\cite{B&R}, Prop. 2.5.22, we obtain that $x_n \rightarrow x$ 
$\sigma$-weakly. 
Moreover, an application of the Hahn-Banach theorem shows that $\m_0 \subset \m_{\phi}$ is a dense subset in the norm topology (see~\cite{B&R} Cor.2.5.23).
This proves the first claim of the lemma. 

Also we deduce from~\cite{B&R}, Prop. 2.5.22,
that an element $A \in \m_0$ is analytic for $\tau_t$, see~\cite{B&R}, Def. 2.5.20. This implies that there exists a strip 
$$ I_{\lambda} = \{ z\ \ | \  |Im z | < \lambda \} \subset \mathbb{C},$$
and a function $f_A : I_{\lambda} \rightarrow \m$, such that 
\begin{equation} 
f_A(t) = \tau_t(A) = \Delta^{it}A\Delta^{-it},
\end{equation}
for $t \in \mathbb{R}$.
In this case, we have $I_{\lambda} = \mathbb{C}$. 
Analyticy for $A \in \m_0$ also implies that $f_A$ 
is an analytic function on 
$\mathbb{C}$ using~\cite{B&R}, Prop. 2.5.21. Hence, we have for $z \in \mathbb{C}$,
\begin{equation}
f_A(z) = \Delta^{iz}A\Delta^{-iz}, 
\end{equation}
since an analytic function is determined by its restriction to the real axis. This proves that $\Delta^{1/2}A\Delta^{-1/2} \in \m$ for $A\in \m_0$. 

It is left to show that $\Delta^{1/2}A\Delta^{-1/2} \in \m_{\phi}$ for $A \in \m_0$. The proof of~\cite{B&R}, Prop. 2.5.22 shows that for $r \in \mathbb{R}, x \in \m_{\phi}$, 
we have 
\begin{equation}
\tau_r(x_n) = \sqrt{n / \pi} \int dt \ e^{-n(t-r)^2} \tau_t(x),
\end{equation}
where the integral on the right hand side is an analytic function on $\mathbb{C}$. Hence, for $z \in \mathbb{C}$, we obtain
\begin{eqnarray}
|\phi(\tau_{z}(x_n)^*\tau_{z}(x_n))| &=& |\frac{\pi}{n} \phi(\int dt\ e^{-n(t-z)^2}\tau_t(x^*) \int ds\ e^{-n(s- z)^2}\tau_s(x))| \nonumber \\
&\leq& \frac{\pi}{n} \int dt\ ds\ e^{-n(t-z)^2} e^{-n(s-z)^2} |\phi(\tau_t(x^*)\tau_s(x))| \nonumber \\
&<& \infty,
\label{finite!}
\end{eqnarray}
by a similar computation as above. Applying \eqref{finite!} to $z = -i/2$, we obtain 
\begin{equation}
\phi(\tau_{1/2}(x_n)^*\tau_{1/2}(x_n)) < \infty.
\end{equation}
Hence, for $A \in \m_0$, we have  
\begin{equation}
\Delta^{1/2}A\Delta^{-1/2} \in \m_{\phi},
\end{equation}

Now we will prove the third claim. 
Applying the construction above to the von Neumann algebra $Z(\m)$ we obtain a dense subset $Z(\m)_0$ of analytic elements, and trivially 
$Z(\m)_0 \subset \m_0$. 
Recall
\begin{equation}
\mathfrak{D}_{\phi} = Span\{x \in \m_+ \ \ | \ \ \phi(x) < \infty \}.
\end{equation}
Now ~\cite{P&T}, Prop. 3.3 states that $\mathfrak{D}_{\phi}$ is a two-sided ideal for $\m_0$, hence 
\begin{eqnarray}
Z(\m)_0\mathfrak{D}_{\phi} &\subset& \mathfrak{D}_{\phi} \nonumber \\
\mathfrak{D}_{\phi}Z(\m)_0 &\subset& \mathfrak{D}_{\phi}.
\end{eqnarray}
It follows from~\cite{P&T}, Thm. 3.6, that $Z(\m)_0$ is contained in the fixed point algebra of the automorphism group $\tau_t$, which means 
\begin{eqnarray}
Z(\m)_0 &\subset& \{ h \in \m \ | \ \ \tau_t(h) = h, t \in \mathbb{R} \} \nonumber \\
&=& \{ h \in \m \ | \ \ \Delta^{it}h\Delta^{-it} = h, t \in \mathbb{R} \}.
\end{eqnarray}
Since $Z(\m)_0$ consists of analytic elements, we may extend the (constant) function $\Delta^{it}h\Delta^{-it}$ to $\mathbb{C}$. 
Hence, for $h \in Z(\m)_0, z \in \mathbb{C}$, we have 
\begin{equation}
\Delta^{iz}h\Delta^{-iz} = h.
\end{equation}
The fact that $Z(\m)_0$ is $\sigma$-dense in $Z(\m)$ finishes the proof.  
\end{bewijs}
\label{lemmaabove}
\end{lemma}

\begin{opmerking}
After this preparation, note that the relative tensor product $\sqtimes_{\n}$ 
has an analogous property  
to the '$S$-balancedness' of the tensor product of bimodules over rings, but with a crucial and interesting modification (see Section \ref{er=tens}). Namely, for $\eta \in \mathcal{H}, 
n \in \n_0, \zeta \in \mathcal{K}$ we have
\begin{equation}
\eta n \boxtimes_{\n} \zeta = \eta \boxtimes_{\n} (\Delta^{1/2}n\Delta^{-1/2})\zeta, \label{equal1}
\end{equation}
and
\begin{equation}
\eta \boxtimes_{\n} n\zeta_1 = \eta (\Delta^{-1/2} n \Delta^{1/2}) \boxtimes_{\n} \zeta. \label{equal2}
\end{equation}
For let $\eta_2 \in \widetilde{\mathcal{H}}, x \in \n_0, \zeta_2 \in \mathcal{K}$; then
\begin{eqnarray}
\Bigl(\eta_1 x \otimes \zeta_1, \eta_2 \otimes \zeta_2\Bigr)_0 &=& <\zeta_1, R_{\eta_1x}^{*\phantom{1x}}R_{\eta_2}^{\phantom{1}}\zeta_2>_{\mathcal{K}} \nonumber \\
&=& <\zeta_1, (\Delta^{1/2}x^*\Delta^{-1/2})R_{\eta_1}^{*\phantom{1}}R_{\eta_2}^{\phantom{1}}\zeta_2>_{\mathcal{K}} \nonumber \\
&=& <(\Delta^{1/2}x\Delta^{-1/2})\zeta_1, R_{\eta_1}^{*\phantom{1n}}R_{\eta_2}^{\phantom{1}}\zeta_2>_{\mathcal{K}} \nonumber \\
&=& \Bigl(\eta_1 \otimes (\Delta^{1/2}x\Delta^{-1/2})\zeta_1, \eta_2 \otimes \zeta_2\Bigr)_0,
\end{eqnarray}
where the second equality holds because of Lemma \ref{c4}.\ref{help}.
This implies that \\$(\eta_1x \otimes \zeta_1) - (\eta_1 \otimes (\Delta^{1/2}x \Delta^{-1/2})\zeta_1)$ belongs to the null space $\mathcal{N}$. 
Since $\n_0 \subset \n_{\psi} \subset \mathfrak{L}^2(\n)$ are all dense, 
equality \eqref{equal1} holds for the completion 
$\mathcal{H} \boxtimes_{\n} \mathcal{K}$. Note that we have used that $\Delta^{1/2}x\Delta^{-1/2}$ is an element of $\n$ for $x \in \n_0$ (see Lemma 
\ref{lemmaabove}), so that the expression $(\Delta^{1/2}x\Delta^{-1/2})\zeta_1$ is defined. The proof of \eqref{equal2} follows by a similar 
argument.
\label{problemsolved!}
\end{opmerking}
\label{reltens}
\end{subsection}

\begin{subsection}{Morita theory with use of bicategories}
As in the case of rings, one would like to express the notion of Morita equivalence of von Neumann algebras in terms of bicategories. First, let us define the representation category of a von Neumann algebra. 

\begin{definitie}
Let $\m$ be a von Neumann algebra. Then $Rep(\m)$ denotes the category of normal unital $^*$-representations on Hilbert spaces as objects, 
and bounded linear intertwiners as arrows.
\label{defmrep}
\end{definitie}

Analogously to the case of rings, one would expect that Morita equivalence of 
von Neumann algebras is defined by categorical equivalence of their representation 
categories. However, we use another definition. 
\begin{definitie}
Two von Neumann algebras $\m,\n$ are called \textdef{Morita equivalent}, $\m \meq \n$, if a 
correspondence $\m \rightarrow \mathcal{H} \leftarrow \n$ exists, where the representation of $\m$ on $\mathcal{H}$ is faithful and for which 
\begin{equation}
\m^{\prime} \simeq \n^{op},
\label{faith}
\end{equation}
holds.
\label{defmoreq}
\end{definitie}
Note that faithfulness of $\m$  directly implies faithfulness for $\n$, using \eqref{faith}. Using somewhat different notation, 
Rieffel (\!~\cite{Rieff}, Thm 8.5) proves equivalence between the definition 
above and equivalence of the representation categories, where the equivalence 
is implemented by a normal $^*$-functor\footnote{In this paper, we will not follow Rieffel's method of proof. His proof involves so-called normal 
$\n$-rigged $\m$-modules, 
and tensor products of these modules. Rieffel proves an Eilenberg-Watts like theorem,(\!~\cite{Rieff}, Thm. 5.5),
 which states that all functors of $Rep(\n)$ to $Rep(\m)$ are equivalent to taking the
 tensor product with such a normal $\n$-rigged $\m$-module. However, these modules are not Hilbert spaces, which is what we would like considering our definition of 
the tensor product and the bicategory of von Neumann algebras.}. See also~\cite{sau}. 
Also, compare Lemma \ref{crucial} to justify this choice of definition.  

Now, we will show that the collection of von Neumann algebras forms a 
 bicategory. This result was already stated in~\cite{land} without proof.

\begin{prop}
For any two von Neumann algebra's $\m$,$\n$, let $(\m,\n)$ be the category of correspondences as objects, and bounded linear bimodule maps as arrows.
Then the collection of all von Neumann algebras as objects and correspondences as arrows forms a bicategory $\mathsf{[W^*]}$. The composition functor 
$(\m,\n) \times (\n,\p) \rightarrow (\m,\p)$ is given by the relative tensor product $\sqtimes_{\n}$, and the unit arrow in $(\m,\m)$ is given by 
$I_{\m} = \m \rightarrow \mathfrak{L}^2(\m) \leftarrow \m$.
\begin{bewijs}
One by one, the properties of a bicategory as stated in section \ref{bic} will be checked.
\begin{itemize}
\item{
The object space $\mathsf{[W^*]}_0$ consists of all von Neumann algebras.}
\item{
For any two von Neumann algebras $\m, \n$, there is a category $(\m,\n)$. $(\m,\n)$ has correspondences $\m \rightarrow \mathcal{H} \leftarrow \n$ as objects 
and bounded linear bimodule maps as arrows.}
\item{
For each triple $(\m,\n,\p)$ of von Neumann algebras there is a composition functor, given by $\sqtimes_{\n}$, defined with an arbitrary normal faithful semifinite weight $\psi$.
On objects $\sqtimes_{\n}$ is defined by: 
\begin{eqnarray} 
 \sqtimes_{\n} : (\m,\n) \times (\n,\p) &\longrightarrow& (\m,\p) \nonumber \\
 _{\m}\mathcal{H}_{\n} \times _{\n}\mathcal{K}_{\p} &\longmapsto& _{\m}(\mathcal{H} \sqtimes_{\n} \mathcal{K} )_{\p},
\end{eqnarray} 
for $_{\m}\mathcal{H}_{\n} \in (\m,\n)_0$ and $ _{\n}\mathcal{K}_{\p} \in (\n,\p)_0$. 
On arrows, $\sqtimes_{\n}$ acts as follows:
Let $f :\  _{\m}{\mathcal{H}_1}_{\n} \rightarrow  _{\m}\!{\mathcal{H}_2}_{\n}$ be an arrow in $(\m,\n)$ and let $g : \ _{\n}{\mathcal{K}_1}_{\p} 
\rightarrow_{\n}\!{\mathcal{K}_2}_{\p}$ be an arrow in $(\n,\p)$. Then, for \hbox{$\sum_i(\eta_i \sqtimes_{\n} \zeta_i)$} $\in \mathcal{H}_1 \sqtimes_{\n}
 \mathcal{K}_1$ we define 
\begin{eqnarray}
f \sqtimes_{\n} g\  :\  _{\m}\!(\mathcal{H}_1 \sqtimes_{\n} \mathcal{K}_1)_{\p} &\longrightarrow& _{\m}(\mathcal{H}_2 \sqtimes_{\n}\mathcal{K}_2)_{\p} 
\nonumber \\  
\sum_i(\eta_i\sqtimes_{\n} \zeta_i) &\longmapsto& \sum_i(f(\eta_i) \sqtimes_{\n} g(\zeta_i)).
\end{eqnarray}
Since both $f$ and $g$ intertwine the $\n$-action (and therefore the \\$\Delta^{1/2}\n_0\Delta^{-1/2}$-action by Lemma \ref{lemmaabove}),
 $f \sqtimes_{\n} g$ is well-defined. 

It has already been shown that on objects, 
the image of $\sqtimes_{\n}$ lies in $(\m,\p)$ (see Section \ref{reltens}). 
On arrows, one has to show that, for each pair of arrows $(f \times g) \in \bigl((\m,\n) \times (\n,\p)\bigr)_1$, the image 
$f \sqtimes_{\n} g$ is a bounded bimodule map (in the sense that it commutes with the left representation of $\m$ and the right representation of 
$\p$). Moreover, it has to be shown that $\sqtimes_{\n}$ satisfies the properties of a functor.
These proofs are left to the reader since they are completely similar to the proofs in the case of the tensor product for rings\footnote{While proving these properties 
in the case of rings, we do not use the $S$-balancedness of the tensor product, which is different in the case of von Neumann algebras (see Remark \ref{problemsolved!}). 
Hence we can apply the same arguments for von Neumann algebras.}, see Section \ref{er=tens}.
}

\item{For each object $\m$ of $\mathsf{[W^*]}_0$, the identity arrow $I_{\m}$ of $(\m,\m)$ is given by $\m \rightarrow \mathfrak{L}^2(\m) \leftarrow \m$, 
the identity correspondence, as discussed above.
}

\item{For each quadruple $(\m,\n,\p,\q)$ of von Neumann algebras, 
we need a natural isomorphism $\beta$ between the functors \hbox{$
((- \sqtimes_{\n} - ) \sqtimes_{\p}- ) $} and \\ \hbox{$ (- \sqtimes_{\n} (- \sqtimes_{\p} -))$}, each from the category  
$(\m,\n) \times (\n,\p) \times (\p,\q)$ to the category $(\m,\q)$. To each object $(\mathcal{H},\mathcal{K},\mathcal{L})$ in 
$(\m,\n) \times (\n,\p) \times (\p,\q)$, let $\beta$ assign an arrow 
$$ (\mathcal{H} \sqtimes_{\n} \mathcal{K} ) \sqtimes_{\p} \mathcal{L} \rightarrow \mathcal{H} \sqtimes_{\n} (\mathcal{K} \sqtimes_{\p}\mathcal{L}),$$ 
in $(\m,\q)$, where 
$\beta$ is defined as follows on the elements:
$$( \eta \sqtimes_{\n} \zeta) \sqtimes_{\p} \xi \mapsto \eta \sqtimes_{\n} (\zeta \sqtimes_{\q} \xi).$$
Now, given an arrow  $$(f,g,h) : (\mathcal{H},\mathcal{K},\mathcal{L}) \rightarrow (\mathcal{H}^{\prime},\mathcal{K}^{\prime},\mathcal{L}^{\prime})$$
 in \hbox{$(\m,\n) \times (\n,\p) \times (\p,\q)$}, the diagram 
\begin{equation}
\xymatrix{
(\mathcal{H} \sqtimes_{\n} \mathcal{K}) \sqtimes_{\p} \mathcal{L} \ar[d]^{(f \sqtimes_{\n} g) \sqtimes_{\p} h}
\ar[rrr]^{\beta(\mathcal{H},\mathcal{K},\mathcal{L})} &&&  \mathcal{H} \sqtimes_{\n} (\mathcal{K} \sqtimes_{\p}\mathcal{L}) 
\ar[d]^{f \sqtimes_{\n} (g \sqtimes_{\p} h)} \\
(\mathcal{H}^{\prime} \sqtimes_{\n} \mathcal{K}^{\prime} ) \sqtimes_{\p} \mathcal{L}^{\prime}
\ar[rrr]^{\beta(\mathcal{H}^{\prime},\mathcal{K}^{\prime},\mathcal{L}^{\prime})} &&&  \mathcal{H}^{\prime} \sqtimes_{\n} (\mathcal{K}^{\prime} 
\sqtimes_{\p} \mathcal{L}^{\prime})} 
\end{equation}
commutes by definition of the relative tensor product on the arrows. Again, observe the similarity to the case of rings.
}

\item{For each pair of von Neumann algebras $(\m,\n)$, we need the left identity $L_{(\m,\n)}$: this should be a natural isomorphism between the functor $$\sqtimes_{\m} \circ 
\bigl[I_{\m} \times Id_{(\m,\n)}\bigr]$$ and the canonical 
functor from $1 \times (\m,\n)$ to $(\m,\n)$. 
 To each object in $ 1 \times (\m,\n)$, let $L_{(\m,\n)}$ assign an arrow 
$$ \mathfrak{L}^2(\m) \sqtimes_{\m} \mathcal{K} \rightarrow \ _{\m}\!\mathcal{K}_{\n},$$ in $(\m,\n)$. 
Let $x \in \m_0, \zeta \in \mathcal{K}$. Then 
$L_{(\m,\n)}( 1 \times _{\m}\!\mathcal{K}_{\n})$ is defined on the dense subspace 
$\m_0 \boxtimes_{\m} \mathcal{K}$ by
\begin{eqnarray}
L_{(\m,\n)}( 1 \times _{\m}\!\mathcal{K}_{\n})\ \  :\ \  
\m_0 \sqtimes_{\m} \mathcal{K}  &\rightarrow& 
_{\m}\mathcal{K}_{\n} \nonumber \\
\Lambda(x) \sqtimes_{\m} \zeta &\mapsto& x\zeta. 
\label{cont1}
\end{eqnarray}

We will show that the map \eqref{cont1} is continuous, so 
that we may extend it to $\mathfrak{L}^2(\m) \sqtimes_{\m} \mathcal{K}$. 
Consider 
\begin{equation}
|| \Lambda(x) \sqtimes_{\m} \zeta ||^2 = <\zeta , R_{\Lambda(x)}^{*\phantom{gF}} R_{\Lambda(x)}^{\phantom{PO}} \zeta >_{\mathcal{K}}. 
\end{equation}
Observe that in this case, the operator $R_{\Lambda(x)} : \mathfrak{L}^2(\m) \rightarrow \mathfrak{L}^2(\m)$ is given by 
\begin{equation}
R_{\Lambda(x)} J\Lambda(y^*) = \pi_r(y)\Lambda(x) = Jy^*J\Lambda(x),
\label{2xhetzelfde}
\end{equation}
on the dense subspace $\m_0$, 
by definition of the right representation of $\m$ on $\mathfrak{L}^2(\m)$. 

Recall that $\Lambda : \m_{\psi} \rightarrow \mathfrak{L}^2(\m)$ denotes the inclusion map. It satisfies 
\begin{eqnarray}
x\Lambda(y) &=& \Lambda(xy) \nonumber \\
\Delta^{1/2}\Lambda(y) &=& \Lambda(\Delta^{1/2} y \Delta^{-1/2}), 
\end{eqnarray}
for $x \in \m, y \in \m_{\psi}$. See~\cite{K&R} or~\cite{P&T}.

The following claims lead to continuity of the map \eqref{cont1}. 
\begin{enumerate}
\item{For $x,z \in \m_0$, we have $R_{\Lambda(x)} \Lambda(z) = \Lambda(xz)$. \\
Proof: We have 
\begin{eqnarray}
R_{\Lambda(x)} J\Lambda(y^*) &=& R_{\Lambda(x)} \Delta^{1/2}S\Lambda(y^*) \nonumber \\ &=& R_{\Lambda(x)} \Delta^{1/2} \Lambda(y) \nonumber \\
&=& R_{\Lambda(x)} \Lambda(\Delta^{1/2}y\Delta^{-1/2}),
\end{eqnarray}
hence, using \eqref{2xhetzelfde} and writing $y = \Delta^{-1/2}z\Delta^{1/2}$, we obtain
\begin{eqnarray}
R_{\Lambda(x)} \Lambda(z) &=& J\Delta^{1/2}z^*\Delta^{-1/2}J\Lambda(x) \nonumber \\
&=& Sz^*S\Lambda(x) = Sz^*\Lambda(x^*) \nonumber \\
&=& S\Lambda(z^*x^*) = \Lambda(xz). 
\label{simpel}
\end{eqnarray}}
\item{For $x \in \m_0$, the operator $R^{*}_xR_x^{\phantom{O}}$ equals $x^*x$ as operators on $\mathfrak{L}^2(\m)$.\\
Proof: Using \eqref{simpel}, we obtain 
\begin{eqnarray}
(\Lambda(y), R^{*}_xR_x^{\phantom{O}}\Lambda(z)) &=& (R_x^{\phantom{c}}\Lambda(y), R_x^{\phantom{c}}\Lambda(z)) \nonumber \\
&=& (\Lambda(xy),\Lambda(xz)) \nonumber \\
&=& \psi(y^*x^*xz) \nonumber \\
&=& (\Lambda(y), x^*x \Lambda(z)),
\end{eqnarray}
for $x,y,z \in \m_0$. }
\end{enumerate}

Combining the statements above we obtain 
\begin{equation}
|| \Lambda(x) \sqtimes_{\m} \zeta ||^2 = <\zeta , R_{\Lambda(x)}^{*\phantom{gF}} R_{\Lambda(x)}^{\phantom{8PO}} \zeta >_{\mathcal{K}} = ||x\zeta||_{\mathcal{K}}^2.
\label{isometrie!}
\end{equation}
Hence $\Lambda(x) \sqtimes \zeta \mapsto x\zeta$ is a continuous map. 
It is left to show that $L_{(\m,\n)}$ is unitary and hence a natural isomorphism. We know that it is isometric, from 
\eqref{isometrie!}. It is sufficient to show that the image of $\m_0 \sqtimes \mathcal{K}$ is dense in $\mathcal{K}$. Generally, $\mathbb{I}$ is not an element of 
$\mathfrak{L}^2(\m)$. But, since $\m_0$ is dense in $\mathfrak{L}^2(\m)$, we may find a net $\mathcal{E}_{\alpha}$ in $\m_0$, converging to $\mathbb{I}$, even in norm. Let
$\zeta \in \mathcal{K}$. Then 
\begin{equation}
||\mathcal{E}_{\alpha}\zeta - \mathbb{I}\zeta||_{\mathcal{K}} \leq ||\mathcal{E}_{\alpha} - \mathbb{I}||_{\m}\cdot  ||\zeta||_{\mathcal{K}} \rightarrow 0,
\end{equation}
where we used the fact that a representation is norm decreasing. 
Hence, $\m_0\mathcal{K}$ is dense in $\mathcal{K}$. 

In contrast to the case of rings, the right identity $R_{(\m,\n)}$ is defined differently from the left identity. The right identity should be a natural isomorphism between the functor 
$$\sqtimes_{\n} \circ \bigl[ Id_{(\m,\n)} \times I_{\n} \bigr] $$ and the canonical functor from $(\m,\n) \times 1$ to $(\m,\n)$.
To each object in $(\m,\n) \times 1$ let $R_{(\m,\n)}$ assign an arrow 
$$ \mathcal{H} \sqtimes_{\n} \mathfrak{L}^2(\n) \rightarrow  \ _{\m}\!\mathcal{H}_{\n}, $$ in $(\m,\n)$. 
Let $\eta \in \widetilde{\mathcal{H}}, y \in \n_0$. Then 
$R_{(\m,\n)}( _{\m}\!\mathcal{H}_{\n} \times 1)$ is defined on the dense 
subspace $\widetilde{\mathcal{H}} \sqtimes_{\n} \n_0$ by
\begin{eqnarray}
R_{(\m,\n)}( _{\m}\!\mathcal{H}_{\n} \times 1)\ \ : \ \ 
\widetilde{\mathcal{H}} \sqtimes_{\n} \n_0 &\rightarrow&  _{\m}\!\mathcal{H}_{\n} \nonumber \\
\eta \sqtimes_{\n} \Lambda(y) &\mapsto& \eta(\Delta^{-1/2}y\Delta^{1/2}).
\label{cont2}
\end{eqnarray}
Note that $\Delta^{-1/2}y\Delta^{1/2} \in \n$, for $y \in \n_0$, so that $\eta(\Delta^{-1/2}y\Delta^{1/2})$ is defined by definition of the right representation of $\n$ on $\mathcal{H}$.
We will show next that the map \eqref{cont2} is continuous, so that we may extend it to $\mathcal{H} \sqtimes_{\n} \mathfrak{L}^2(\n)$.
Consider
\begin{equation}
||\eta \sqtimes_{\n} \Lambda(y)||^2 = (\Lambda(y), R_{\eta}^{*}R_{\eta}^{\phantom{P}}\Lambda(y))_{\mathfrak{L}^2(\n)}.
\end{equation}
In this case, the operator $R_{\eta} : \mathfrak{L}^2(\n) \rightarrow \mathcal{H}$ is defined by 
\begin{equation}
R_{\eta}J\Lambda(y^*) = \pi_r(y)\eta = \eta y,
\end{equation}
on the dense subspace $\n_0$. 

We make use of the following fact, proven by Connes (\!~\cite{Conn2} Lemma 4). For $\xi \in \widetilde{\mathcal{H}}$, we have 
\begin{equation}
\psi(R_{\xi}^*R_{\xi}^{\phantom{P}}) = ||\xi||^2_{\mathcal{H}}.\end{equation} 
Hence, using Lemma \ref{5e}, we obtain 
\begin{eqnarray}
||\eta \sqtimes_{\n} \Lambda(y)||^2 &=&(\Lambda(y), R_{\eta}^{*}R_{\eta}^{\phantom{P}}\Lambda(y))_{\mathfrak{L}^2(\n)} \nonumber \\
&=& \psi(y^*R_{\eta}^*R_{\eta}^{\phantom{I}}y) \nonumber \\
&=& \psi(R_{\eta(\Delta^{-1/2}y\Delta^{1/2})}^{*\phantom{boel}}R_{\eta(\Delta^{-1/2}y\Delta^{1/2})}^{\phantom{boelbo}}) \nonumber \\
&=& ||\eta(\Delta^{-1/2}y\Delta^{1/2})||^2_{\mathcal{H}}.
\label{isometrie22!!}
\end{eqnarray}
Hence $\eta \sqtimes_{\n} \Lambda(y) \mapsto \eta(\Delta^{1/2}y\Delta^{-1/2})$ is a continuous map.
It is left to show that $R_{(\m,\n)}$ is unitary and hence a natural isomorphism. As above, \eqref{isometrie22!!} shows that 
it is isometric, hence it is sufficient to show that the image of $\widetilde{\mathcal{H}} \sqtimes \n_0$ is dense in $\mathcal{H}$. As before, we have a net 
$\mathcal{E}_{\alpha}$ in $\n_0$ converging to \nolinebreak $\mathbb{I}$. Consider the net $\Delta^{1/2}\mathcal{E}_{\alpha}\Delta^{-1/2}$. This net is contained in 
$\n_{\psi}$, following from Lemma \ref{lemmaabove}. By the inclusion $\n_0 \subset \n_{\psi} \subset \mathfrak{L}^2(\n)$ and the continuity just proven, 
we have 
\begin{equation}
\Delta^{1/2}\mathcal{E}_{\alpha}\Delta^{-1/2} \sqtimes_{\n} \eta \mapsto \eta\mathcal{E}_{\alpha},
\end{equation} 
for $\eta \in \widetilde{\mathcal{H}}$. 
The right hand side converges to $\eta$ in norm. The fact that $\widetilde{\mathcal{H}} \subset \mathcal{H}$ is a dense subspace finishes the proof. 
}

\item{
We need to prove associativity coherence. Let $\m,\n,\p,\q,\mathfrak{R}$ be 
von Neumann algebras, and $_{\m}\mathcal{J}_{\n}\ ,\ _{\n}\mathcal{H}_{\p}\ ,\
 _{\p}\mathcal{K}_{\q}\ ,\ _{\q}\mathcal{L}_{\mathfrak{R}}$ be associated correspondences. The following diagram
\begin{equation}
\xymatrix{
((\mathcal{J} \sqtimes_{\n} \mathcal{H}) \sqtimes_{\p} \mathcal{K}) \sqtimes_{\q} \mathcal{L} 
\ar[rrr]^{\beta(\mathcal{J},\mathcal{H},\mathcal{K}) \sqtimes_{\mathfrak{R}}Id}
\ar[d]_{\beta(\mathcal{J}\sqtimes_{\n}\mathcal{H},\mathcal{K},\mathcal{L})} 
&&&(\mathcal{J} \sqtimes_{\n} (\mathcal{H} \sqtimes_{\p} \mathcal{K})) \sqtimes_{\q} \mathcal{L} \ar[d]^{\beta(\mathcal{J},\mathcal{H}\sqtimes_{\p}\mathcal{K},\mathcal{L})} \\
(\mathcal{J} \sqtimes_{\n} \mathcal{H}) \sqtimes_{\p} (\mathcal{K} \sqtimes_{\q} \mathcal{L}) 
\ar[ddrrr]_{\beta(\mathcal{J},\mathcal{H},\mathcal{K}\sqtimes_{\q} \mathcal{L})\ \ \ \ \ \ }
&&& \mathcal{J} \sqtimes_{\n} ((\mathcal{H} \sqtimes_{\p} \mathcal{K} ) 
\sqtimes_{\q} \mathcal{L}) \ar[dd]^{Id \sqtimes_{\n} \beta(\mathcal{H},\mathcal{K},\mathcal{L})} \\ 
&&&& \\
 &&& \mathcal{J} \sqtimes_{\n} (\mathcal{H} \sqtimes_{\p} (\mathcal{K} \sqtimes_{\q} \mathcal{L})), &&\\
}
\end{equation}
commutes, since $\beta$ is defined elementwise.
}
\item{
We need to prove identity coherence. Let $\m,\n,\p$ be von Neumann algebras and $_{\m}\mathcal{H}_{\n}, _{\n}\mathcal{K}_{\p}$ associated bimodules. 
It will be shown that the following diagram commutes:
\begin{equation}
\xymatrix{
(\mathcal{H} \sqtimes_{\n} \mathfrak{L}^2(\n)) \sqtimes_{\n} \mathcal{K} \ar[rr]^{\beta(\mathcal{H},\mathfrak{L}^2(\n),\mathcal{K})}, 
\ar[dr]_{R(\mathcal{H},\mathfrak{L}^2(\n)) \sqtimes_{\n} Id\ \ \ \ \ } 
&& \mathcal{H} \sqtimes_{\n} ( \mathfrak{L}^2(\n) \sqtimes_{\n} \mathcal{K})\ar[dl]^{\ \ \ \ \ Id \sqtimes_{\n} L(\mathfrak{L}^2(\n),\mathcal{K})} \\
& \mathcal{H} \sqtimes_{\n} \mathcal{K}. & \\
}
\label{diag2}
\end{equation}
Since the upper arrow $\beta(\mathcal{H},\mathfrak{L}^2(\n),\mathcal{K})$ is the associativity isomorphism, identity coherence boils down to \begin{equation}
(\eta (\Delta^{-1/2} n \Delta^{1/2}) \sqtimes_{\n} \zeta) \cong (\eta \sqtimes_{\n} n\zeta),
\label{noversq}
\end{equation}
for $\eta \in \mathcal{H}$, $n \in \mathfrak{L}^2(\n)$, $\zeta \in \mathcal{K}$, which follows immediately from Remark \ref{problemsolved!}.
}
\end{itemize}
\end{bewijs}
\label{vna=bic}
\end{prop}  
Note that the bicategory of von Neumann algebras constructed as above, depends on the choice of the faithful, semifinite, normal weight for each von Neumann algebra. The definition of the composition functor (i.e. the relative tensor product), depends on the choice of the weight. However, different weights lead to unitary equivalent relative tensor products. Hence the bicategory depends on the weight only up to isomorphism.

The following theorem is an analogue of the Morita theorem for rings. 
\begin{stelling}
Two von Neumann algebras are isomorphic in the bicategory $\mathsf{[W^*]}$ iff they are Morita equivalent. In formula:
\begin{equation}
\m \biso \n \Longleftrightarrow \m \meq \n.
\end{equation}

\begin{bewijs}
First, let us reformulate the first statement of the theorem. Recall Definition \ref{defisobic}. 
The property that two von Neumann algebras $\m,\n$ are isomorphic in the bicategory means that 
there must be an arrow in $(\m,\n)$, i.e., a correspondence $\m \rightarrow \mathcal{H} \leftarrow \n$ that is invertible up to isomorphism. In other words, 
there exists an arrow $\n \rightarrow \mathcal{H}^{-1} \leftarrow \m$ in $(\n,\m)$ such that 
\begin{eqnarray}
\m \rightarrow \mathcal{H} \sqtimes_{\n} \mathcal{H}^{-1} \leftarrow \m &\simeq& \m \rightarrow \mathfrak{L}^2(\m) \leftarrow \m 
\text{\ \ in\ \ } (\m,\m) 
\label{inv1}\\
\n \rightarrow \mathcal{H}^{-1} \sqtimes_{\m} \mathcal{H} \leftarrow \n &\simeq& \n \rightarrow \mathfrak{L}^2(\n) \leftarrow \n 
\text{\ \ in\ \ } (\n,\n). \label{inv2}
\end{eqnarray}
``$\Leftarrow:$'' Let $\m,\n$ be Morita equivalent. So we have a 
correspondence \\ $\m \rightarrow \mathcal{H} \leftarrow \n$, where the representation of $\m$ on $\mathcal{H}$ is faithful, and $\m^{\prime} \simeq \n^{op}$.
>From $\m \rightarrow \mathcal{H} \leftarrow \n$, we can define a correspondence 
 $\n \rightarrow \overline{\mathcal{H}} \leftarrow \m$ by 
\begin{equation}
n\bar{\eta} m := m^*\eta n^*, \hspace{10pt} \text{for \ } n \in \n, m \in \m, 
\eta \in \mathcal{H},
\end{equation}
where $\overline{\mathcal{H}}$ is $\mathcal{H}$ as a set, with the addition operator of $\mathcal{H}$ and conjugate scalar multiplication and inner product.
Sauvageot~\cite{sau}, Prop 3.1 proves that the relative tensor product 
$\mathcal{H} \sqtimes_{\n} \overline{\mathcal{H}}$ is in standard form (using explicitly that \\
$\m^{\prime} \simeq \n^{op}$), i.e.,
\begin{equation}
\m \rightarrow \mathcal{H} \sqtimes_{\n} \overline{\mathcal{H}} \leftarrow \m \hspace{10pt}\simeq \hspace{10pt} \m \rightarrow \mathfrak{L}^2(\m) \leftarrow \m.
\label{onehand}
\end{equation}
Applying the same reasoning to $\n \rightarrow \overline{\mathcal{H}} \leftarrow \m$,  we obtain $\m \rightarrow \overline{\overline{\mathcal{H}}} \leftarrow \n$ and 
clearly, we have $\overline{\overline{\mathcal{H}}} = \mathcal{H}$. Then 
\begin{eqnarray}
\n \rightarrow \overline{\mathcal{H}} \sqtimes_{\m} \overline{\overline{\mathcal{H}}} \leftarrow \n &=&
\n \rightarrow \overline{\mathcal{H}} \sqtimes_{\m} \mathcal{H} \leftarrow \n \nonumber \\
&\cong& \n \rightarrow \mathfrak{L}^2(\n) \leftarrow \n.
\label{otherhand}
\end{eqnarray} 
Together, \eqref{onehand} and \eqref{otherhand} prove that $_{\m}\mathcal{H}_{\n}$ is invertible, its inverse being $_{\n}\overline{\mathcal{H}}_{\m}$.

``$\Rightarrow:$'' Suppose we have an invertible correspondence $\m \rightarrow \mathcal{H} \leftarrow \n$. We need to show that $\m \simeq (\n^{op})^{\prime}$, and that 
the representation of $\m$ is faithful. By definition of a correspondence, we have 
\begin{equation}
\m \subseteq (\n^{op})^{\prime},
\label{inclus}
\end{equation}
so, considering the representation of $\m$ on $\mathcal{H} \sqtimes_{\n} \mathcal{H}^{-1}$, one has 
\begin{equation}
\m \sqtimes id_{\mathcal{H}^{-1}} \subseteq (\n^{op})^{\prime} \sqtimes id_{\mathcal{H}^{-1}}.
\end{equation}
Now we will use a result from Sauvageot~\cite{sau}, Prop. 3.3, who shows that for a von Neumann algebra $\p$ and representations 
$\mathcal{K}_1 \leftarrow \p$ and $\p \rightarrow \mathcal{K}_2$, one has
\begin{eqnarray}
\bigl[(\p^{op})^{\prime} \sqtimes_{\p} id_{\mathcal{K}_2} \bigr]^{\prime} &=& id_{\mathcal{K}_1} \sqtimes_{\p} \p^{\prime} \label{sau1} \\
\text{and}  \nonumber \\
(\p^{op})^{\prime} \sqtimes_{\p} id_{\mathcal{K}_2} &=& \bigl[id_{\mathcal{K}_1} \sqtimes_{\p} \p^{\prime} \bigr]^{\prime} \label{sau2}
\end{eqnarray}
in $\mathcal{K}_1 \sqtimes_{\p} \mathcal{K}_2$.
Applying \eqref{sau2} we obtain
\begin{equation}
\m \sqtimes_{\n} id_{\mathcal{H}^{-1}} \subseteq (\n^{op})^{\prime} \sqtimes_{\n} id_{\mathcal{H}^{-1}} = \bigl[id_{\mathcal{H}} \sqtimes_{\n} \n^{\prime} \bigr]^{\prime},
\end{equation}
hence, using \eqref{sau1},
\begin{equation}
id_{\mathcal{H}} \sqtimes_{\n} \n^{\prime} \subseteq \bigl[\m \sqtimes_{\n} id_{\mathcal{H}^{-1}}\bigr]^{\prime} = id_{\mathcal{H}} \sqtimes_{n} \m^{op}.
\label{final}
\end{equation}
Equation \eqref{final} implies that $\n^{\prime} \subseteq \m^{op}$: Suppose $\n^{\prime} \nsubseteq \m^{op}$. Then an element $n\in \n^{\prime} $ exists, such that $\forall m \in \m^{op}$:
$$\exists\ \ \  \eta_1 \sqtimes_{\n} \eta_2 \in \mathcal{H} \sqtimes_{\n} \mathcal{H}^{-1},$$
such that 
\begin{equation}
(\eta_1 \sqtimes_{\n} n \eta_2 ) \neq (\eta_1 \sqtimes_{\n} m \eta_2).
\end{equation}
This would violate \eqref{final}.
Hence
\begin{equation}
(\n^{op})^{\prime} \subseteq \m,
\end{equation}
which, with \eqref{inclus}, proves $\m \simeq (\n^{op})^{\prime}$. 
It remains to be shown that the representation of $\m$ on $\mathcal{H}$ is 
faithful. However, this follows immediately from \eqref{inv1} and the fact that the standard representation of $\m$ on $\mathfrak{L}^2(\m)$ is faithful.
\end{bewijs}
\label{main}
\end{stelling}

It is possible to restate Theorem \ref{main} above in terms of representation 
categories. In the light of the remarks after Definition \ref{defmoreq}, the proof is immediate. However, we will prove the theorem directly, c.f. the proof of 
Theorem \ref{Bicat}.

\begin{stelling}
Two von Neumann algebras are isomorphic in the bicategory $\mathsf{[W^*]}$ iff their representation categories are equivalent, where the equivalence is 
implemented by a normal $^*$-functor. 
\begin{bewijs}
``$\Rightarrow :$'' Given the invertible correspondence $_{\m}\mathcal{H}_{\n}$ we will construct a functor $F : Rep(\m) \rightarrow Rep(\n)$ and a functor $G : Rep(\n) \rightarrow Rep(\m)$ such that $F \circ G \cong id_{Rep(\n)}$ and $G \circ F \cong id_{Rep(\m)}$. 
On objects, define
\begin{equation}
F_0(\mathcal{K}) := {\mathcal{H}}^{-1} \sqtimes_{\m} \mathcal{K}, \hspace{10pt} \mathcal{K} \in Rep(\m)_0.
\end{equation}
On arrows, define
\begin{equation}
F_1(f) := id_{{\mathcal{H}}^{-1}} \sqtimes_{\m} f, \hspace{10pt}  f \in Rep(\m)_1.
\end{equation}
Since ${\mathcal{H}}^{-1}$ is a $\n-\m$ correspondence, 
the image of $F_0$ lies in $Rep(\n)_0$. We will check that $F$ satisfies the 
properties of a functor (see \eqref{propfunct1} and \eqref{propfunct2}). Let $\eta \in \mathcal{K} \in Rep(\m)_0 , \zeta \in \mathcal{H}^{-1}$. 
Then one has
\begin{eqnarray}
F_1(id_{\mathcal{K}})(\zeta \sqtimes_{\m} \eta) &=& \bigl(id_{{\mathcal{H}}^{-1}} \sqtimes_{\m} id_{\mathcal{K}}\bigr) (\zeta \sqtimes_{\m} \eta) \nonumber \\
&=&\bigl( id_{{\mathcal{H}}^{-1}}(\zeta) \sqtimes_{\m} id_{\mathcal{K}}(\eta)\bigr) = (\zeta \sqtimes_{\m} \eta) \nonumber \\
&=& id_{({\mathcal{H}^{-1}} \sqtimes_{\m} \mathcal{K})} (\zeta \sqtimes_{\m} \eta) \nonumber \\
&=& id_{F_1(\mathcal{K})}(\zeta \sqtimes_{\m} \eta),
\end{eqnarray}
and
\begin{eqnarray} 
F_1(h_1 \circ h_2)(\zeta \sqtimes_{\m} \eta) &=& id_{{\mathcal{H}}^{-1}}(\zeta) \sqtimes_{\m} (h_1 \circ h_2)(\eta)\nonumber \\
&=& F_1(h_1) \bigl[(\zeta \sqtimes_{\m} h_2(\eta))\bigr] \nonumber \\
&=& F_1(h_1) \circ F_1(h_2) \bigl[(\zeta \sqtimes_{\m} \eta)\bigr];
\end{eqnarray}
for $h_1,h_2 \in Rep(\m)_1$, where the last equation holds whenever $h_1 \circ h_2$ is defined in $Rep(\m)_1$.

In the same way, one constructs a functor $G : Rep(\n) \rightarrow Rep(\m),$ by putting $G_0(\mathcal{L}) := {\mathcal{H}} \sqtimes_{\n} \mathcal{L} $ on objects, and $G_1(g) := id_{{\mathcal{H}}} \sqtimes_{\n} g $ on arrows. 

Now we will show that natural equivalences $(F \circ G) \cong id_{Rep(\n)}$ and \\ $(G \circ F) \cong id_{Rep(\m)}$ exist. Using equation \eqref{inv1} and the 
fact that the relative tensor product is associative up to isomorphism, we obtain
\begin{eqnarray}
(G \circ F)_0(\mathcal{K}) &=& G_0({\mathcal{H}^{-1}} \sqtimes_{\m} \mathcal{K}) = \mathcal{H} \sqtimes_{\n} ({\mathcal{H}}^{-1} \sqtimes_{\m} {\mathcal{K}}) 
\nonumber \\ &\cong& (\mathcal{H} \sqtimes_{\n} {\mathcal{H}}^{-1}) \sqtimes_{\m} {\mathcal{K}} \simeq  \mathcal{L}^2(\m) \sqtimes_{\m} 
\mathcal{K} \nonumber \\  &\cong& 
\mathcal{K},
\label{congruent}
\end{eqnarray}
for $\mathcal{K} \in Rep(\m)_0$. 
Let $f$ be an arrow $\mathcal{K}_1 \rightarrow \mathcal{K}_2$ in $Rep(\m)_1$. 
Via the isomorphism in \eqref{congruent}, the diagram
\begin{equation}
\xymatrix{
 \mathcal{K}_1 \ar[d]^f \ar[rr] && (G \circ F)(\mathcal{K}_1) \ar[d]^{id_{\mathcal{H}} \sqtimes_{\n} (id_{{\mathcal{H}}^{-1}} \sqtimes_{\m} f)} \\
\mathcal{K}_2 \ar[rr]&& \mathcal{H} \sqtimes_{\n} ({\mathcal{H}}^{-1} \sqtimes_{\m} \mathcal{K}_2),
}
\end{equation}
commutes, so $(G \circ F) \cong id_{Rep(\m)}$. A similar computation shows that \\ $(F \circ G) \cong id_{Rep(\n)}$. 

It is left to show that $F,G$ are normal $^*$-functors. 
 We will show that $F$ is a $^*$-functor, i.e. ${F_1(f)}^* = F_1(f^*)$. Let $f : \mathcal{K}_1 \rightarrow \mathcal{K}_2$. 
Then for $\eta_1,\eta_2 \in \widetilde{\mathcal{H}^{-1}}, \zeta_1 \in \mathcal{K}_2, \zeta_2 \in \mathcal{K}_1$, one has  
\begin{eqnarray}
\bigl(F_1(f)^* (\eta_1 \otimes \zeta_1), \eta_2 \otimes \zeta_2\bigr)_0 &=& \bigl(\eta_1 \otimes \zeta_1 , F_1(f)(\eta_2 \otimes \zeta_2)\bigr)_0 \nonumber \\
&=& \bigl(\eta_1 \otimes \zeta_1, \eta_2 \otimes f(\zeta_2)\bigr)_0 \nonumber \\
&=& <\zeta_1, R_{\eta_1}^*R_{\eta_2}^{*\phantom{1}} f(\zeta_2)>_{\mathcal{K}_2} \nonumber \\
&=& <\zeta_1, f(R_{\eta_1}^*R_{\eta_2}^{*\phantom{1}} \zeta_2)>_{\mathcal{K}_2} \hspace{15pt} (\star) \nonumber  \\
&=& <f^*(\zeta_1) , R_{\eta_1}^*R_{\eta_2}^{*\phantom{1}}\zeta_2>_{\mathcal{K}_2} \nonumber \\
&=& \bigl(\eta_1 \otimes f^*(\zeta_1), \eta_2 \otimes \zeta_2\bigr)_0 \nonumber \\
&=& \bigl(F_1(f^*)(\eta_1 \otimes \zeta_1), \eta_2 \otimes \zeta_2\bigr)_0. 
\label{above}
\end{eqnarray}
The equality $(\star)$ holds, since $f$ intertwines the $\n$-action.
Note that, since $f$ maps into $\mathcal{K}_2$, the inner product $(\ ,\ )_0$ is defined on $\widetilde{\mathcal{H}^{-1}} \otimes_{\mathbb{C}} \mathcal{K}_2$. 
Equation \eqref{above} implies that $F_1(f)^* = F_1(f^*)$ on the closure $\mathcal{H}^{-1} \sqtimes \mathcal{K}_1$. 
Since $F$ is a 
category equivalence (up to isomorphism), it is 
automatically full and faithful. It then follows from~\cite{Rieff}, Prop 7.3 
that $F$ is normal. A similar argument holds for $G$. 

``$\Leftarrow:$'' 
Suppose $Rep(\n) \simeq Rep(\m)$, where the categorical equivalence is implemented by a $^*$-functor.  Let $F : Rep(\m) \rightarrow Rep(\n)$. 
Consider $F(\mathfrak{L}^2(\m))$, which has a left $\n$-action by definition, and a right $\m$-action through $F_1$.
Sauvageot has stated an Eilenberg-Watts like theorem,~\cite{sau}, Prop 5.3, also cf. Rieffel~\cite{Rieff}, Prop 5.4 and Thm 5.5:
\begin{stelling}[Sauvageot]
An equivalence functor $F : Rep(\m) \rightarrow Rep(\n)$ is characterized by the $\n$-$\m$ bimodule $F(\mathfrak{L}^2(\m))$ in the following way. 
On objects,
\begin{equation}
F_0(\mathcal{K}) \cong F(\mathfrak{L}^2(\m)) \sqtimes_{\m} \mathcal{K},
\end{equation}
for $\mathcal{K} \in Rep(\m)_0$.
On arrows, for $f : \mathcal{K}_1 \rightarrow \mathcal{K}_2$,
\begin{equation}
F_1(f) \cong id_{\mathfrak{L}^2(\m)} \sqtimes f.
\end{equation}
Hence
\begin{equation}
F(-) \cong  F(\mathfrak{L}^2(\m)) \sqtimes_{\m} (-) .
\end{equation}
\end{stelling}
Similarly, for $G : Rep(\n) \rightarrow Rep(\m)$, we have,
\begin{equation}
G(-) \cong G(\mathfrak{L}^2(\n)) \sqtimes_{\n} (-) .
\end{equation}

Using the fact that $G \circ F \cong id_{Rep(\m)}$ and $F \circ G \cong id_{Rep(\n)}$ and that composition of natural equivalences does provide a natural equivalence, we now find
\begin{equation}
F(\mathfrak{L}^2(\m)) \sqtimes_{\n} G(\mathfrak{L}^2(\n)) \cong F\bigl[G(\mathfrak{L}^2(\n))\bigr] \cong \mathfrak{L}^2(\n),
\end{equation}
and
\begin{equation}
G(\mathfrak{L}^2(\n)) \sqtimes_{\m} F(\mathfrak{L}^2(\m)) \cong G\bigl[F(\mathfrak{L}^2(\m))\bigr] \cong \mathfrak{L}^2(\m),
\end{equation}
as bimodules.
This shows that $F(\mathfrak{L}^2(\m))$ is an invertible $\n$-$\m$ bimodule: its inverse (up to isomorphism) is $G(\mathfrak{L}^2(\n))$.
\end{bewijs}
\label{last}
\end{stelling}
\end{subsection}
\end{section}

%% file: bibl.tex
\addcontentsline{toc}{section}{References}
\bibliographystyle{alpha}
  